\documentclass{amsart}

\usepackage{a4wide,amssymb,amsrefs,calligra,graphicx,mathrsfs,txfonts}

\usepackage[2cell,arrow,curve,matrix,web]{xy}
\UseHalfTwocells
\UseTwocells
\newdir{ >}{{}*!/-10pt/\dir{>}}

\def\calli#1{\textup{\!\textcalligra{#1}\,}}

\def\A{{\mathscr A}}
\def\a{{\mathbf a}}
\def\B{{\mathbf B}}
\def\b{{\mathbf b}}
\def\aB{{\mathscr B}}

\def\E{{\mathrm E}}
\def\F{{\mathbb F}}
\def\G{{\mathbb G}}

\def\H{{\mathscr H}}
\def\L{{\mathrm L}}
\def\LL{{\mathbb L}}
\def\P{{\mathrm P}}
\def\R{{\mathrm R}}
\def\RR{{\mathbb R}}
\def\Z{{\mathbb Z}}
\def\ZZ{{\mathbb Z}}

\def\Topt{{\mathbf{Top}^*}}
\def\Sp{{\mathbf{Spec}}}

\def\op{^{\mathrm{op}}}
\def\ho{_\simeq}
\def\1{^{-1}}
\def\comp{{\scriptstyle\Box}}

\def\Pair{{\calli{Pair}}}
\def\Mod{{\mathbf{Mod}}}

\DeclareMathOperator{\Aut}{Aut}
\DeclareMathOperator{\Ext}{Ext}
\DeclareMathOperator{\Hom}{Hom}
\DeclareMathOperator{\im}{im}
\DeclareMathOperator{\coker}{coker}
\DeclareMathOperator{\Sq}{Sq}
\DeclareMathOperator{\Tr}{Tr}

\def\xto#1{\xrightarrow[]{#1}}
\def\xot#1{\xleftarrow[]{#1}}
\def\set#1{\left\{#1\right\}}

\def\hog#1{\left\llbracket#1\right\rrbracket}

\def\brk#1{\left\langle#1\right\rangle}
\def\alignbox#1{\begin{aligned}#1\end{aligned}}
\def\smat#1{{\left(\begin{smallmatrix}#1\end{smallmatrix}\right)}}
\def\restr#1#2{\left.#1\right|_{#2}}

\let\ge\geqslant
\let\le\leqslant
\let\d\partial

\let\ox\otimes
\let\ot\leftarrow
\let\then\Rightarrow
\let\onto\twoheadrightarrow
\let\into\hookrightarrow

\makeatletter

\let\c@equation=\c@subsection

\def\thmhead#1#2#3{%
  (\thmnumber{
    \@upn{#2}})
   \thmname{#1}%
   \thmnote{ {\the\thm@notefont(#3)}}}

\makeatother

\def\subsection#1{\refstepcounter{subsection}\medskip\noindent{\textbf{(\thesubsection)
\ #1\unskip. }}\ignorespaces}

\def\subsubsection#1{\refstepcounter{subsubsection}\smallskip\noindent{\textbf{(\thesubsubsection)}
\ \textit{#1\unskip. }}\ignorespaces}

\def\thesubsection{\thesection.\arabic{subsection}}
\def\thesubsubsection{\thesection.\arabic{subsection}.\arabic{subsubsection}}

\def\theodef#1{\newtheorem{#1}[subsection]{#1}}

\theodef{Lemma}
\theodef{Proposition}
\theodef{Theorem}
\theodef{Corollary}
\newtheorem{Theocite}[subsection]{Theorem \cite{Baues}}

\theoremstyle{definition}

\theodef{Definition}
\theodef{Remark}
\theodef{Example}
\theodef{Examples}
\theodef{Construction}

\begin{document}

\title{Computation of the $E_3$-term of the Adams spectral sequence}

\author{Hans Joachim Baues}
\address{
Max-Planck-Institut f\"ur Mathematik\\
Vivatsgasse 7\\
D-53111 Bonn\\
Germany}
\email{baues@mpim-bonn.mpg.de}

\author{Mamuka Jibladze}
\address{
Razmadze Mathematical Institute\\
M.~Alexidze st. 1\\
Tbilisi 0193\\
Georgia}
\email{jib@rmi.acnet.ge}

\maketitle

The algebra $\aB$ of secondary cohomology operations is a pair algebra with
$\Sigma$-structure which as a Hopf algebra was explicitly computed in
\cite{Baues}. In particular the multiplication map $A$ of $\aB$ was
determined by an algorithm. In this paper we introduce algebraically the
secondary $\Ext$-groups $\Ext_\aB$ and we prove that the E$_3$-term of the
Adams spectral sequence (computing stable maps in $\set{Y,X}^*_p$) is given
by
$$
\E_3(Y,X)=\Ext_\aB(\H X,\H Y).
$$
Here $\H X$ is the secondary cohomology of the spectrum $X$ which is the
$\aB$-module $\G^\Sigma$ if $X$ is the sphere spectrum $S^0$. This leads to
an algorithm for the computation of the group
$$
\E_3(S^0,S^0)=\Ext_\aB(\G^\Sigma,\G^\Sigma)
$$
which is a new explicit approximation of stable homotopy groups of spheres
improving the Adams approximation
$$
\E_2(S^0,S^0)=\Ext_\A(\F,\F).
$$
An implementation of our algorithm computed $\E_3(S^0,S^0)$ by now up
to degree 40. In this range our results confirm the known results in the
literature, see for example the book of Ravenel \cite{Ravenel}.

\section{Modules over pair algebras}

We here recall from \cite{Baues} the notion of pair modules, pair
algebras, and pair modules over a pair algebra $B$. The category
$B$-$\Mod$ of pair modules over $B$ is an additive track category in which
we consider secondary resolutions as defined in \cite{BJ5}. Using such
secondary resolutions we shall obtain the secondary derived functors
$\Ext_B$ in section \ref{secodif}.

Let $k$ be a commutative ring with unit and let $\Mod$ be the category of
$k$-modules and $k$-linear maps. This is a symmetric monoidal category via
the tensor product $A\!\ox\!B$ over $k$ of $k$-modules $A$, $B$. A \emph{pair}
of modules is a morphism
\begin{equation}\label{pair}
X=\left(X_1\xto\d X_0\right)
\end{equation}
in $\Mod$. We write $\pi_0(X)=\ker\d$ and $\pi_1(X)=\coker\d$. A
\emph{morphism} $f:X\to Y$ of pairs is a commutative diagram
$$
\xymatrix
{
X_1\ar[r]^{f_1}\ar[d]_\d&Y_1\ar[d]^\d\\
X_0\ar[r]^{f_0}&Y_0.
}
$$
Evidently pairs with these morphisms form a category $\Pair(\Mod)$ and
one has functors
$$
\pi_0, \pi_1 : \Pair(\Mod)\to\Mod.
$$
A pair morphism is called a
\emph{weak equivalence} if it induces isomorphisms on $\pi_0$ and $\pi_1$.

Clearly a pair in $\Mod$ coincides with a chain complex concentrated in
degrees 0 and 1. For two pairs $X$ and $Y$ the tensor product of the
complexes corresponding to them is concentrated in degrees in 0, 1 and 2
and is given by
$$
X_1\!\ox\!Y_1\xto{\d_1}X_1\!\ox\!Y_0\oplus
X_0\!\ox\!Y_1\xto{\d_0}X_0\!\ox\!Y_0
$$
with $\d_0=(\d\ox1,1\ox\d)$ and $\d_1=(-1\ox\d,\d\ox1)$. Truncating $X\ox Y$
we get the pair
$$
X\bar\otimes Y=
\left((X\bar\otimes Y)_1=\coker(\d_1)\xto\d X_0\ox Y_0=(X\bar\otimes Y)_0\right)
$$
with $\d$ induced by $\d_0$.

\begin{Remark}
Note that the full embedding of the category of pairs into the category of
chain complexes induced by the above identification has a left adjoint
$\Tr$ given by truncation: for a chain complex
$$
C=\left(...\to C_2\xto{\d_1}C_1\xto{\d_0}C_0\xto{\d_{-1}}C_{-1}\to...\right),
$$
one has
$$
\Tr(C)=\left(\coker(\d_1)\xto{\bar{\d_0}}C_0\right),
$$
with $\bar{\d_0}$ induced by $\d_0$. Then clearly one has
$$
X\bar\otimes Y=\Tr(X\ox Y).
$$
Using the fact that $\Tr$ is a reflection onto a full subcategory, one
easily checks that the category $\Pair(\Mod)$ together with the tensor
product $\bar\otimes$ and unit $k=(0\to k)$ is a symmetric monoidal category,
and $\Tr$ is a monoidal functor.
\end{Remark}

We define the tensor product $A\ox B$ of two graded modules in the usual
way, i.~e. by
$$
(A\ox B)^n=\bigoplus_{i+j=n}A^i\ox B^j.
$$

A \emph{(graded) pair module} is a graded object of $\Pair(\Mod)$, i.~e. a
sequence $X^n=(\d:X_1^n\to X_0^n)$ of pairs in $\Mod$. We identify such a
pair module $X$ with the underlying morphism $\d$ of degree 0 between
graded modules
$$
X=\left(X_1\xto\d X_0\right).
$$
Now the tensor product $X\bar\otimes Y$ of graded pair modules $X$, $Y$ is defined by 
\begin{equation}\label{grpr}
(X\bar\otimes Y)^n=\bigoplus_{i+j=n}X^i\bar\otimes Y^j.
\end{equation}
This defines a monoidal structure on the category of graded pair modules.
Morphisms in this category are of degree 0.

For two morphisms $f,g:X\to Y$ between graded pair modules, a
\emph{homotopy} $H:f\then g$ is a morphism $H:X_0\to Y_1$ of degree 0 as in
the diagram
\begin{equation}\label{homot}
\alignbox{
\xymatrix
{
X_1\ar@<.5ex>[r]^{f_1}\ar@<-.5ex>[r]_{g_1}\ar[d]_\d&Y_1\ar[d]^\d\\
X_0\ar@<.5ex>[r]^{f_0}\ar@<-.5ex>[r]_{g_0}\ar[ur]|H&Y_0,
}}
\end{equation}
satisfying $f_0-g_0=\d H$ and $f_1-g_1=H\d$.

A \emph{pair algebra} $B$ is a monoid in the monoidal category of graded pair
modules, with multiplication
$$
\mu:B\bar\otimes B\to B.
$$
We assume that $B$ is concentrated in nonnegative degrees, that is $B^n=0$
for $n<0$.

A \emph{left $B$-module} is a graded pair module $M$ together with a left
action
$$
\mu:B\bar\otimes M\to M
$$
of the monoid $B$ on $M$.

More explicitly pair algebras and modules over them can be described as
follows.

\begin{Definition}
A \emph{pair algebra} $B$ is a graded pair
$$
\d:B_1\to B_0
$$
in $\Mod$ with $B_1^n=B_0^n=0$ for $n<0$ such that $B_0$ is a graded
algebra in $\Mod$, $B_1$ is a graded $B_0$-$B_0$-bimodule, and $\d$ is a
bimodule homomorphism. Moreover for $x,y\in B_1$ the equality
$$
\d(x)y=x\d(y)
$$
holds in $B_1$.
\end{Definition}

It is easy to see that there results an exact sequence of graded
$B_0$-$B_0$-bimodules
$$
0\to\pi_1B\to B_1\xto\d B_0\to\pi_0B\to0
$$
where in fact $\pi_0B$ is a $k$-algebra, $\pi_1B$ is a
$\pi_0B$-$\pi_0B$-bimodule, and $B_0\to\pi_0(B)$ is a homomorphism of
algebras.

\begin{Definition}\label{bmod}
A \emph{(left) module} over a pair algebra $B$ is a graded pair
$M=(\d:M_1\to M_0)$ in $\Mod$ such that $M_1$ and $M_0$ are left
$B_0$-modules and $\d$ is $B_0$-linear. Moreover a $B_0$-linear map
$$
\bar\mu:B_1\!\otimes_{B_0}\!M_0\to M_1
$$
is given fitting in the commutative diagram
$$
\xymatrix{
B_1\otimes_{B_0}M_1\ar[r]^{1\ox\d}\ar[d]_\mu&B_1\otimes_{B_0}M_0\ar[dl]^{\bar\mu}\ar[d]^\mu\\
M_1\ar[r]_\d&M_0,
}
$$
where $\mu(b\ox m)=\d(b)m$ for $b\in B_1$ and $m\in M_1\cup M_0$.

For an indeterminate element $x$ of degree $n=|x|$ let $B[x]$ denote the
$B$-module with $B[x]_i$ consisting of expressions $bx$ with $b\in B_i$,
$i=0,1$, with $bx$ having degree $|b|+n$, and structure maps given by
$\d(bx)=\d(b)x$, $\mu(b'\ox bx)=(b'b)x$ and $\bar\mu(b'\ox bx)=(b'b)x$.

A \emph{free} $B$-module is a direct sum of several copies of modules of
the form $B[x]$, with $x\in I$ for some set $I$ of indeterminates of
possibly different degrees. It will be denoted
$$
B[I]=\bigoplus_{x\in I}B[x].
$$

For a left $B$-module $M$ one has the exact sequence of $B_0$-modules
$$
0\to\pi_1M\to M_1\to M_0\to\pi_0M\to0
$$
where $\pi_0M$ and $\pi_1M$ are actually $\pi_0B$-modules.

Let $B$-$\Mod$ be the category of left modules over the pair algebra $B$.
Morphisms $f=(f_0,f_1):M\to N$ are pair morphisms which are
$B$-equivariant, that is,$f_0$ and $f_1$ are $B_0$-equivariant and 
compatible with $\bar\mu$ above, i.~e. the diagram
$$
\xymatrix{
B_1\ox_{B_0}M_0\ar[r]^-{\bar\mu}\ar[d]_{1\ox f_0}&M_1\ar[d]^{f_1}\\
B_1\ox_{B_0}N_0\ar[r]^-{\bar\mu}&N_1
}
$$
commutes.

For two such maps $f,g:M\to N$ a track $H:f\then g$ is a degree zero map
\begin{equation}\label{track}
H:M_0\to N_1
\end{equation}
satisfying $f_0-g_0=\d H$ and $f_1-g_1=H\d$ such that $H$ is
$B_0$-equivariant. For tracks $H:f\then g$, $K:g\then h$ their composition
$K\comp H:f\then h$ is $K+H$.
\end{Definition}

\begin{Proposition}
For a pair algebra $B$, the category $B$-$\Mod$ with the above track
structure is a well-defined additive track category.
\end{Proposition}

\begin{proof}
For a morphism $f=(f_0,f_1):M\to N$ between $B$-modules, one has
$$
\Aut(f)=\set{H\in\Hom_{B_0}(M_0,N_1)\ |\ \d H=f_0-f_0,H\d=f_1-f_1}\cong\Hom_{\pi_0B}(\pi_0M,\pi_1N).
$$
Since this group is abelian, by \cite{Baues&JibladzeI} we know that $B$-$\Mod$
is a linear track extension of its homotopy category by the bifunctor $D$
with $D(M,N)=\Hom_{\pi_0B}(\pi_0M,\pi_1N)$. It thus remains to show that the
homotopy category is additive and the bifunctor $D$ is biadditive.

By definition the set of morphisms $[M,N]$ between objects $M$, $N$ in the
homotopy category is given by the exact sequence of abelian groups
$$
\Hom_{B_0}(M_0,N_1)\to\Hom_B(M,N)\onto[M,N].
$$
This makes evident the abelian group structure on the hom-sets $[M,N]$. Bilinearity
of composition follows from consideration of the commutative diagram
$$
\xymatrix{
\Hom_{B_0}(M_0,N_1)\!\ox\!\Hom_B(N,P)\oplus
\Hom_B(M,N)\!\ox\!\Hom_{B_0}(N_0,P_1)\ar[d]\ar[r]^-\mu
&\Hom_{B_0}(M_0,P_1)\ar[d]\\
\Hom_B(M,N)\ox\Hom_B(N,P)\ar[r]\ar@{->>}[d]
&\Hom_B(M,P)\ar@{->>}[d]\\
[M,N]\ox[N,P]\ar@{-->}[r]
&[M,P]
}
$$
with exact columns, where $\mu(H\!\ox\!g+f\!\ox\!K)=g_1H+Kf_0$. It also shows
that the functor $B$-$\Mod\to B$-$\Mod\ho$ is linear. Since this functor
is the identity on objects, it follows that the homotopy category is additive.

Now note that both functors $\pi_0$, $\pi_1$ factor to define functors on
$B$-$\Mod\ho$. Since these functors are evidently additive, it follows that
$D=\Hom_{\pi_0B}(\pi_0,\pi_1)$ is a biadditive bifunctor.
\end{proof}

\begin{Lemma}\label{freehom}
If $M$ is a free $B$-module, then the canonical map
$$
[M,N]\to\Hom_{\pi_0B}(\pi_0M,\pi_0N)
$$
is an isomorphism for any $B$-module $N$. 
\end{Lemma}

\begin{proof}
Let $(g_i)_{i\in I}$ be a free generating set for $M$. Given a
$\pi_0(B)$-equivariant homomorphism $f:\pi_0M\to\pi_0N$, define its lifting
$\tilde f$ to $M$ by specifying $\tilde f(g_i)=n_i$, with $n_i$ chosen
arbitrarily from the class $f([g_i])=[n_i]$.

To show monomorphicity, given $f:M\to N$ such that $\pi_0f=0$, this means
that $\im f_0\subset\im\d$, so we can choose $H(g_i)\in N_1$ in such a way
that $\d H(g_i)=f_0(g_i)$. This then extends uniquely to a $B_0$-module
homomorphism $H:M_0\to N_1$ with $\d H=f_0$; moreover any element of $M_1$ is a linear
combination of elements of the form $b_1g_i$ with $b_1\in B_1$, and for
these one has $H\d(b_1g_i)=H(\d(b_1)g_i)=\d(b_1)H(g_i)$. But
$f_1(b_1g_i)=b_1f_0(g_i)=b_1\d H(g_i)=\d(b_1)H(g_i)$ too, so $H\d=f_1$.
This shows that $f$ is nullhomotopic.
\end{proof}

\section{$\Sigma$-structure}

\begin{Definition}
The \emph{suspension} $\Sigma X$ of a graded object $X=(X^n)_{n\in\Z}$ is
given by degree shift, $(\Sigma X)^n=X^{n-1}$.
\end{Definition}

Let $\Sigma:X\to\Sigma X$ be the map of degree 1 given by the identity. If
$X$ is a left $A$-module over the graded algebra $A$ then $\Sigma X$ is a
left $A$-module via
\begin{equation}\label{suspact}
a\cdot\Sigma x=(-1)^{|a|}\Sigma(a\cdot x)
\end{equation}
for $a\in A$, $x\in X$. On the other hand if $X$ is a right $A$-module then
$(\Sigma x)\cdot a=\Sigma(x\cdot a)$ yields the right $A$-module structure
on $\Sigma X$.

\begin{Definition}\label{sigma}
A \emph{$\Sigma$-module} is a graded pair module $X=(\d:X_1\to X_0)$
together with an isomorphism
$$
\sigma:\pi_1X\cong\Sigma\pi_0X
$$
of graded $k$-modules. We then call $\sigma$ a \emph{$\Sigma$-structure} of
$X$. A $\Sigma$-map between $\Sigma$-modules is a map $f$ between pair
modules such that $\sigma(\pi_1f)=\Sigma(\pi_0f)\sigma$. If $X$ is a pair
algebra then a $\Sigma$-structure is an isomorphism of
$\pi_0X$-$\pi_0X$-bimodules. If $X$ is a left module over a pair algebra $B$
then a $\Sigma$-structure of $X$ is an isomorphism $\sigma$ of left
$\pi_0B$-modules. Let
$$
(B\textrm{-}\Mod)^\Sigma\subset B\textrm{-}\Mod
$$
be the track category of $B$-modules with $\Sigma$-structure and
$\Sigma$-maps.
\end{Definition}

\begin{Lemma}
Suspension of a $B$-module $M$ has by \eqref{suspact} the structure of a
$B$-module and $\Sigma M$ has a $\Sigma$-structure if $M$ has one.
\end{Lemma}

\begin{proof}
Given $\sigma:\pi_1M\cong\Sigma\pi_0M$ one defines a $\Sigma$-structure on
$\Sigma M$ via
$$
\pi_1\Sigma
M=\Sigma\pi_1M\xto{\Sigma\sigma}\Sigma\Sigma\pi_0M=\Sigma\pi_0\Sigma M.
$$
\end{proof}

Hence we get suspension functors between track categories
$$
\xymatrix{
B\textrm{-}\Mod\ar[r]^\Sigma&B\textrm{-}\Mod\\
(B\textrm{-}\Mod)^\Sigma\ar[u]\ar[r]^\Sigma&(B\textrm{-}\Mod)^\Sigma.\ar[u]
}
$$

\begin{Lemma}\label{add}
The track category $(B\mathrm{-}\Mod)^\Sigma$ is $\LL$-additive in the sense of
\cite{BJ5}, with $\LL=\Sigma\1$, or as well $\RR$-additive,
with $\RR=\Sigma$.
\end{Lemma}

\begin{proof}
The statement of the lemma means that the bifunctor
$$
D(M,N)=\Aut(0_{M,N})
$$
is either left- or right-representable, i.~e. there is an endofunctor
$\LL$, respectively $\RR$ of $(B$-$\Mod)^\Sigma$ and a binatural
isomorphism $D(M,N)\cong[\LL M,N]$, resp. $D(M,N)\cong[M,\RR N]$.

Now by \eqref{track}, a track in $\Aut(0_{M,N})$ is a $B_0$-module
homomorphism $H:M_0\to N_1$ with $\d H=H\d=0$; hence
$$
D(M,N)\cong\Hom_{\pi_0B}(\pi_0M,\pi_1N)\cong\Hom_{\pi_0B}(\pi_0\Sigma\1
M,\pi_0N)\cong\Hom_{\pi_0B}(\pi_0M,\pi_0\Sigma N).
$$
\end{proof}

\begin{Lemma}
If $B$ is a pair algebra with $\Sigma$-structure then each free $B$-module
has a $\Sigma$-structure.
\end{Lemma}

\begin{proof}
This is clear from the description of free modules in \ref{bmod}.
\end{proof}

\section{The secondary differential over pair algebras}\label{secodif}

For a pair algebra $B$ with a $\Sigma$-structure, for a $\Sigma$-module
$M$ over $B$, and a module $N$ over $B$ we now define the \emph{secondary differential}
$$
d_{(2)}:\Ext^n_{\pi_0B}(\pi_0M,\pi_0N)\to\Ext^{n+2}_{\pi_0B}(\pi_0M,\pi_1N).
$$
Here $d_{(2)}=d_{(2)}(M,N)$ depends on the $B$-modules $M$ and $N$ and is
natural in $M$ and $N$ with respect to maps in $(B\mathrm{-}\Mod)^\Sigma$. For the
definition of $d_{(2)}$ we consider secondary chain complexes and secondary
resolutions. In \cite{BJ5} such a construction was performed in the
generality of an arbitrary $\LL$-additive track category. We will first
present the construction of $d_{(2)}$ for the track category of pair
modules and then will indicate how this construction is a particular case
of the more general situation discussed in \cite{BJ5}.

\begin{Definition}\label{secs}
For a pair algebra $B$, a \emph{secondary chain complex} $M_\bullet$ in $B$-$\Mod$
is given by a diagram of the form
$$
\xymatrix@!{
...\ar[r]
&M_{n+2,1}\ar[r]^{d_{n+1,1}}\ar[d]|\hole^>(.75){\d_{n+2}}
&M_{n+1,1}\ar[r]^{d_{n,1}}\ar[d]|\hole^>(.75){\d_{n+1}}
&M_{n,1}\ar[r]^{d_{n-1,1}}\ar[d]|\hole^>(.75){\d_n}
&M_{n-1,1}\ar[r]\ar[d]|\hole^>(.75){\d_{n-1}}
&...\\
...\ar[r]\ar[urr]^<(.3){H_{n+1}}
&M_{n+2,0}\ar[r]_{d_{n+1,0}}\ar[urr]^<(.3){H_n}
&M_{n+1,0}\ar[r]_{d_{n,0}}\ar[urr]^<(.3){H_{n-1}}
&M_{n,0}\ar[r]_{d_{n-1,0}}\ar[urr]
&M_{n-1,0}\ar[r]
&...\\
}
$$
where each $M_n=(\partial_n:M_{n,1}\to M_{n,0})$ is a $B$-module, each
$d_n=(d_{n,0},d_{n,1})$ is a morphism in $B$-$\Mod$, each $H_n$ is
$B_0$-linear and moreover the identities
\begin{align*}
d_{n,0}d_{n+1,0}&=\d_nH_n\\
d_{n,1}d_{n+1,1}&=H_n\d_{n+2}\\
\intertext{and}
H_nd_{n+2,0}&=d_{n,1}H_{n+1}
\end{align*}
hold for all $n\in\Z$. We thus see that in this case a secondary complex is
the same as a graded version of a \emph{multicomplex} (see e.~g. \cite{Meyer}) with only two
nonzero rows.

One then defines the \emph{total complex} Tot$(M_\bullet)$ of the form
$$
...\ot
M_{n-1,0}\oplus M_{n-2,1}
\xot{
\left(
\begin{smallmatrix}
d_{n-1,0}&-\d_{n-1}\\
H_{n-2}&-d_{n-2,1}
\end{smallmatrix}
\right)
}
M_{n,0}\oplus M_{n-1,1}
\xot{
\left(
\begin{smallmatrix}
d_{n,0}&-\d_n\\
H_{n-1}&-d_{n-1,1}
\end{smallmatrix}
\right)
}
M_{n+1,0}\oplus M_{n,1}
\ot
...
$$
Cycles and boundaries in this complex will be called secondary cycles,
resp. secondary boundaries of $M_\bullet$. Thus a secondary $n$-cycle in
$M_\bullet$ is a pair $(c,\gamma)$ with $c\in M_{n,0}$, $\gamma\in
M_{n-1,1}$ such that $d_{n-1,0}c=\d_{n-1}\gamma$,
$H_{n-2}c=d_{n-2,1}\gamma$ and such a cycle is a
boundary iff there exist $b\in M_{n+1,0}$ and $\beta\in M_{n,1}$ with
$c=d_{n,0}b+\d_n\beta$ and $\gamma=H_{n-1}b+d_{n-1,1}\beta$. A secondary
complex $M_\bullet$ is called \emph{exact} if its total complex is, that is, if
secondary cycles are secondary boundaries.
\end{Definition}

Let us now consider a secondary chain complex $M_\bullet$ in $B$-$\Mod$. It
is clear then that
$$
...\to\pi_0M_{n+2}\xto{\pi_0d_{n+1}}\pi_0M_{n+1}\xto{\pi_0d_n}\pi_0M_n\xto{\pi_0d_{n-1}}\pi_0M_{n-1}\to...
\leqno{\pi_0M_\bullet:}
$$
is a chain complex of $\pi_0B$-modules. The next result corresponds to
\cite{BJ5}*{lemma 3.5}.

\begin{Proposition}\label{exact}
Let $M_\bullet$ be a secondary complex consisting of $\Sigma$-modules and
$\Sigma$-maps between them. If $\pi_0(M_\bullet)$ is an exact complex then
$M_\bullet$ is an exact secondary complex. Conversely, if $\pi_0M_\bullet$ is
bounded below then secondary exactness of $M_\bullet$ implies exactness of
$\pi_0M_\bullet$.
\end{Proposition}

\begin{proof} The proof consists in translating the argument from the
analogous general statement in \cite{BJ5} to our setting.
Suppose first that $\pi_0M_\bullet$ is an exact complex, and consider a
secondary cycle $(c,\gamma)\in M_{n,0}\oplus M_{n-1,1}$, i.~e. one has
$d_{n-1,0}c=\d_{n-1}\gamma$ and $H_{n-2}c=d_{n-2,1}\gamma$. Then in
particular $[c]\in\pi_0M_n$ is a cycle, so there exists
$[b]\in\pi_0M_{n+1}$ with $[c]=\pi_0(d_n)[b]$. Take $b\in[b]$, then
$c-d_{n,0}b=\d_n\beta$ for some $\beta\in M_{n+1,1}$. Consider
$\delta=\gamma-H_{n-1}b-d_{n-1,1}\beta$. One has
$\d_{n-1}\delta=\d_{n-1}\gamma-\d_{n-1}H_{n-1}b-\d_{n-1}d_{n-1,1}\beta=d_{n-1,0}c-d_{n-1,0}d_{n,0}b-d_{n-1,0}\d_n\beta=0$,
so that $\delta$ is an element of $\pi_1M_n$. Moreover
$d_{n-2,1}\delta=d_{n-2,1}\gamma-d_{n-2,1}H_{n-1}b-d_{n-2,1}d_{n-1,1}\beta=H_{n-2}c-H_{n-2}d_{n,0}b-H_{n-2}\d_n\beta=0$,
i.~e. $\delta$ is a cycle in $\pi_1M_\bullet$. Since by assumption
$\pi_0M_\bullet$ is exact, taking into account the $\Sigma$-structure
$\pi_1M_\bullet$ is exact too, so that there exists $\psi\in\pi_1M_n$
with $\delta=d_{n-1,1}\psi$. Define $\tilde\beta=\beta+\psi$. Then
$d_{n,0}b+\d_n\tilde\beta=d_{n,0}b+\d_n\beta=c$ since $\psi\in\ker\d_n$.
Moreover
$d_{n-1,1}\tilde\beta=d_{n-1,1}\beta+d_{n-1,1}\psi=d_{n-1,1}\beta+\delta=\gamma-H_{n-1}b$,
which means that $(c,\gamma)$ is the boundary of $(b,\tilde\beta)$. Thus
$M_\bullet$ is an exact secondary complex.

Conversely suppose $M_\bullet$ is exact, and $\pi_0M_\bullet$ bounded
below. Given a cycle $[c]\in\pi_0(M_n)$, represent it by a $c\in M_{n,0}$.
Then $\pi_0d_{n-1}[c]=0$ implies $d_{n-1,0}c\in\im\d_{n-1}$, so there is a
$\gamma\in M_{n-1,1}$ such that $d_{n-1,0}c=\d_{n-1}\gamma$. Consider
$\omega=d_{n-2,1}\gamma-H_{n-2}c$. One has
$\d_{n-2}\omega=\d_{n-2}d_{n-2,1}\gamma-\d_{n-2}H_{n-2}c=d_{n-2,0}\d_{n-1}\gamma-d_{n-2,0}d_{n-1,0}c=0$,
i.~e. $\omega$ is an element of $\pi_1M_{n-2}$. Moreover
$d_{n-3,1}\omega=d_{n-3,1}d_{n-2,1}\gamma-d_{n-3,1}H_{n-2}c=H_{n-3}\d_{n-1}\gamma-H_{n-3}d_{n,0}c=0$,
so $\omega$ is a $n-2$-dimensional cycle in $\pi_1M_\bullet$. Using the
$\Sigma$-structure, this then gives a $n-3$-dimensional cycle in
$\pi_0M_\bullet$. Now since $\pi_0M_\bullet$ is bounded below, we might
assume by induction that it is exact in dimension $n-3$, so that $\omega$
is a boundary. That is, there exists $\alpha\in\pi_1M_{n-1}$ with
$d_{n-2,1}\alpha=\omega$. Define $\tilde\gamma=\gamma-\alpha$; then one
has
$d_{n-2,1}\tilde\gamma=d_{n-2,1}\gamma-d_{n-2,1}\alpha=d_{n-2,1}\gamma-\omega=H_{n-2}c$.
Moreover $\d_{n-1}\tilde\gamma=\d_{n-1}\gamma=d_{n-1,0}c$ since
$\alpha\in\ker(\d){n-1}$. Thus $(c,\tilde\gamma)$ is a secondary cycle,
and by secondary exactness of $M_\bullet$ there exists a pair $(b,\beta)$
with $c=d_{n,0}b+\d_n\beta$. Then $[c]=\pi_0(d_n)[b]$, i.~e. $c$ is a
boundary.
\end{proof}

\begin{Definition}
Let $B$ be a pair algebra with $\Sigma$-structure.
A \emph{secondary resolution} of a $\Sigma$-module $M=(\d:M_1\to M_0)$ over
$B$ is an exact secondary complex $F_\bullet$ in $(B\mathrm{-}\Mod)^\Sigma$ of the form
$$
\xymatrix@!{
...\ar[r]\ar@{}[d]|\cdots
&F_{31}\ar[r]^{d_{21}}\ar[d]|\hole^>(.75){\d_3}
&F_{21}\ar[r]^{d_{11}}\ar[d]|\hole^>(.75){\d_2}
&F_{11}\ar[r]^{d_{01}}\ar[d]|\hole^>(.75){\d_1}
&F_{01}\ar[r]^{\epsilon_1}\ar[d]|\hole^>(.75){\d_0}
&M_1\ar[d]|\hole^>(.9)\d\ar[r]
&0\ar[r]\ar[d]|\hole
&0\ar[r]\ar[d]
&...\ar@{}[d]|\cdots
\\
...\ar[r]\ar[urr]^<(.3){H_2}
&F_{30}\ar[r]_{d_{20}}\ar[urr]^<(.3){H_1}
&F_{20}\ar[r]_{d_{10}}\ar[urr]^<(.3){H_0}
&F_{10}\ar[r]_{d_{00}}\ar[urr]^<(.3){\hat\epsilon}
&F_{00}\ar[r]_{\epsilon_0}\ar[urr]
&M_0\ar[r]\ar[urr]
&0\ar[r]
&0\ar[r]
&...
}
$$
where each $F_n=(\d_n:F_{n1}\to F_{n0})$ is a free $B$-module.
\end{Definition}

It follows from \ref{exact} that for any secondary resolution $F_\bullet$
of a $B$-module $M$ with $\Sigma$-structure, $\pi_0F_\bullet$ will be a
free resolution of the $\pi_0B$-module $\pi_0M$, so that in particular one
has
$$
\Ext^n_{\pi_0B}(\pi_0M,U)\cong H^n\Hom(\pi_0F_\bullet,U)
$$
for all $n$ and any $\pi_0B$-module $U$.

\begin{Definition}\label{secod}
Given a pair algebra $B$ with $\Sigma$-structure, a $\Sigma$-module $M$
over $B$, a module $N$ over $B$ and a secondary resolution $F_\bullet$ of $M$, we define
the \emph{secondary differential}
$$
d_{(2)}:\Ext^n_{\pi_0B}(\pi_0M,\pi_0N)\to\Ext^{n+2}_{\pi_0B}(\pi_0M,\pi_1N)
$$
in the following way. Suppose given a class
$[c]\in\Ext^n_{\pi_0B}(\pi_0M,\pi_0N)$.
First represent it by some element in $\Hom_{\pi_0B}(\pi_0F_n,\pi_0N)$
which is a cocycle, i.~e. its composite with $\pi_0(d_n)$ is 0. By
\ref{freehom} we know that the natural maps
$$
[F_n,N]\to\Hom_{\pi_0B}(\pi_0F,\pi_0N)
$$
are isomorphisms, hence to any such element corresponds a homotopy class
in $[F_n,N]$ which is also a cocycle, i.~e. value of $[d_n,N]$ on it is
zero. Take a representative map $c:F_n\to N$ from this homotopy class. Then
$cd_n$ is nullhomotopic, so we can find a $B_0$-equivariant map $H:F_{n+1,0}\to N_1$
such that in the diagram
$$
\xymatrix@!{
&F_{n+2,1}\ar[r]^{d_{n+1,1}}\ar[d]_{\d_{n+2}}
&F_{n+1,1}\ar[r]^{d_{n,1}}\ar[d]|\hole^>(.75){\d_{n+1}}
&F_{n,1}\ar[r]^{c_1}\ar[d]|\hole^>(.75){\d_n}
&N_1\ar[d]^\d\\
&F_{n+2,0}\ar[r]_{d_{n+1,0}}\ar[urr]^<(.3){H_n}
&F_{n+1,0}\ar[r]_{d_{n,0}}\ar[urr]^<(.3)H
&F_{n,0}\ar[r]_{c_0}
&N_0.
}
$$
one has $c_0d_{n,0}=\d H$, $c_1d_{n,1}=H\d_{n+1}$ and $\d c_1=c_0\d_n$.
Then taking $\Gamma=c_1H_n-Hd_{n+1,0}$ one has $\d\Gamma=0$,
$\Gamma\d_{n+2}=0$, so $\Gamma$ determines a map
$\bar\Gamma:\coker\d_{n+2}\to\ker\d$, i.~e. from $\pi_0F_{n+2}$ to
$\pi_1N$. Moreover $\bar\Gamma\pi_0(d_{n+2})=0$, so it is a cocycle in
$\Hom(\pi_0(F_\bullet),\pi_1N)$ and we define
$$
d_{(2)}[c]=[\bar\Gamma]\in\Ext^{n+2}_{\pi_0B}(\pi_0M,\pi_1N).
$$
\end{Definition}

\begin{Definition}\label{secodef}
Let $M$ and $N$ be $B$-modules with $\Sigma$-structure. Then also all
the $B$-modules $\Sigma^kM$, $\Sigma^kN$ have $\Sigma$-structures and we
get by \ref{secod} the secondary differential
$$
\xymatrix@!C=16em{
\Ext^n_{\pi_0B}(\pi_0M,\pi_0\Sigma^kN)\ar@{=}[d]\ar[r]^{d_{(2)}(M,\Sigma^kN)}&\Ext^{n+2}_{\pi_0B}(\pi_0M,\pi_1\Sigma^kN)\ar@{=}[d]\\
\Ext^n_{\pi_0B}(\pi_0M,\Sigma^k\pi_0N)\ar[r]^d&\Ext^{n+2}_{\pi_0B}(\pi_0M,\Sigma^{k+1}\pi_0N).
}
$$
In case the composite
$$
\Ext^{n-2}_{\pi_0B}(\pi_0M,\Sigma^{k-1}\pi_0N)\xto d\Ext^n_{\pi_0B}(\pi_0M,\Sigma^k\pi_0N)\xto d\Ext^{n+2}_{\pi_0B}(\pi_0M,\Sigma^{k+1}\pi_0N)
$$
vanishes we define the \emph{secondary $\Ext$-groups} to be the quotient
groups
$$
\Ext^n_B(M,N)^k:=\ker d/\im d.
$$
\end{Definition}

\begin{Theorem}\label{sigmacoinc}
For a $\Sigma$-algebra $B$, a $B$-module $M$ with $\Sigma$-structure and
any $B$-module $N$, the secondary differential $d_{(2)}$ in \ref{secod}
coincides with the secondary differential
$$
d_{(2)}:\Ext^n_\a(M,N)\to\Ext^{n+2}_\a(M,N)
$$
from \cite{BJ5}*{Section 4} as constructed for the $\LL$-additive track
category $(B\mathrm{-}\Mod)^\Sigma$ in \ref{add}, relative to the
subcategory $\b$ of free $B$-modules with $\a=\b\ho$.
\end{Theorem}

\begin{proof}
We begin by recalling the appropriate notions from \cite{BJ5}. There
secondary chain complexes $A_\bullet=(A_n,d_n,\delta_n)_{n\in\Z}$ are
defined in arbitrary additive track category $\B$. They consist of objects
$A_n$, morphisms $d_n:A_{n+1}\to A_n$ and tracks
$\delta_n:d_nd_{n+1}\then0_{A_{n+2},A_n}$, $n\in\Z$, such that the
equality of tracks
$$
\delta_nd_{n+2}=d_n\delta_{n+1}
$$
holds for all $n$. For an object $X$, an $X$-valued $n$-cycle in a
secondary chain complex $A_\bullet$ is defined to be a pair $(c,\gamma)$
consisting of a morphism $c:X\to A_n$ and a track
$\gamma:d_{n-1}c\then0_{X,A_{n-1}}$ such that the equality of tracks
$$
\delta_{n-2}c=d_{n-2}\gamma
$$
is satisfied. Such a cycle is called a \emph{boundary} if there exists a
map $b:X\to A_{n+1}$ and a track $\beta:c\then d_nb$ such that the equality
$$
\gamma=\delta_{n-1}b\comp d_{n-1}\beta
$$
holds. A secondary chain complex is called $X$-exact if every $X$-valued
cycle in it is a boundary. Similarly it is called \emph{$\b$-exact}, if it
is $X$-exact for every object $X$ in $\b$, where $\b$ is a track
subcategory of $\B$. A secondary $\b$-resolution of an object $A$ is a
$\b$-exact secondary chain complex $A_\bullet$ with $A_n=0$ for $n<-1$,
$A_{-1}=A$, $A_n\in\b$ for $n\ne-1$; the last differentials will be then
denoted $d_{-1}=\epsilon:A_0\to A$, $\delta_{-1}=\hat\epsilon:\epsilon
d_0\to0_{A_1,A}$ and the pair $(\epsilon,\hat\epsilon)$ will be called
\emph{augmentation} of the resolution. It is clear that any secondary
chain complex $(A_\bullet,d_\bullet,\delta_\bullet)$ in $\B$ gives rise to
a chain complex $(A_\bullet,[d_\bullet])$, in the ordinary sense, in the
homotopy category $\B\ho$ of $\B$. Moreover if $\B$ is $\Sigma$-additive,
i.~e. there exists a functor $\Sigma$ and isomorphisms
$\Aut(0_{X,Y})\cong[\Sigma X,Y]$, natural in $X$, $Y$, then $\b$-exactness
of $(A_\bullet,d_\bullet,\delta_\bullet)$ implies $\b\ho$-exactness of
$(A_\bullet,[d_\bullet])$ in the sense that the chain complex of abelian
groups $[X,(A_\bullet,[d_\bullet])]$ will be exact for each $X\in\b$. In
\cite{BJ5}, the notion of $\b\ho$-relative derived functors has been
developed using such $\b\ho$-resolutions, which we also recall.

For an additive subcategory $\a=\b\ho$ of the homotopy category $\B\ho$,
the $\a$-relative left derived functors $\L^\a_nF$, $n\ge0$, of a
functor $F:\B\ho\to\A$ from $\B\ho$ to an abelian category $\mathscr
A$ are defined by
$$
(\L^\a_nF)A=H_n(F(A_\bullet)),
$$
where $A_\bullet$ is given by any $\a$-resolution of $A$. Similarly,
$\a$-relative right derived functors of a contravariant functor
$F:\B\ho\op\to\A$ are given by
$$
(\R_\a^nF)A=H^n(F(A_\bullet)).
$$
In particular, for the contravariant functor $F=[\_,B]$ we get the $\a$-relative $\Ext$-groups
$$
\Ext^n_\a(A,B):=(\R_\a^n[\_,B])A=H^n([A_\bullet,B])
$$
for any $\a$-exact resolution $A_\bullet$ of $A$.
Similarly, for the contravariant functor $\Aut(0_{\_,B})$ which assigns to
an object $A$ the group $\Aut(0_{A,B})$ of all tracks
$\alpha:0_{A,B}\then0_{A,B}$ from the zero map $A\to*\to B$ to itself, one
gets the groups of $\a$-derived automorphisms
$$
\Aut^n_\a(A,B):=(\R^n_\a\Aut(0_{\_,B}))(A).
$$

It is proved in \cite{BJ5} that under mild conditions (existence of a
subset of $\a$ such that every object of $\a$ is a direct summand of a
direct sum of objects from that subset) every object has an
$\a$-resolution, and that the resulting groups do not depend on the choice
of a resolution. 

We next recall the construction of the secondary differential from
\cite{BJ5}. This is the map of the form
$$
d_{(2)}:\Ext^n_\a(A,B)\to\Aut^n_\a(0_{A,B});
$$
it is constructed from any secondary $\b$-resolution
$(A_\bullet,d_\bullet,\delta_\bullet,\epsilon,\hat\epsilon)$ of the object
$A$. Given an element $[c]\in\Ext^n_\a(A,B)$, one first represents it
by an $n$-cocycle in $[(A_\bullet,[d_\bullet]),B]$, i.~e. by a homotopy
class $[c]\in[A_n,B]$ with $[cd_n]=0$. One then chooses an actual
representative $c:A_n\to B$ of it in $\B$ and a track $\gamma:0\then
cd_n$. It can be shown that the composite track
$\Gamma=c\delta_n\comp\gamma d_{n+1}\in\Aut(0_{A_{n+2},B})$ satisfies
$\Gamma d_{n+1}=0$, so it is an $(n+2)$-cocycle in the cochain complex
$\Aut(0_{(A_\bullet,[d_\bullet]),B})\cong[(\Sigma A_\bullet,[\Sigma
d_\bullet]),B]$, so determines a cohomology class
$d{(2)}([c])=[\Gamma]\in\Ext^{n+2}_\a(\Sigma A,B)$. It is proved in
\cite{BJ5}*{4.2} that the above construction does not indeed depend on
choices.

Now turning to our situation, it is straightforward to verify that a
secondary chain complex in the sense of \cite{BJ5} in the track category
$B$-$\Mod$ is the same as the 2-complex in the sense of \ref{secs}, and
that the two notions of exactness coincide. In particular then the notions
of resolution are also equivalent.

The track subcategory $\b$ of free modules is generated by coproducts from
a single object, so $\b\ho$-resolutions of any $B$-module exist. In fact
it follows from \cite{BJ5}*{2.13} that any $B$-module has a secondary
$\b$-resolution too.

Moreover there are natural isomorphisms
$$
\Aut(0_{M,N})\cong\Hom_{\pi_0B}(\pi_0M,\pi_1N).
$$
Indeed a track from the zero map to itself is a $B_0$-module homomorphism
$H:M_0\to N_1$ with $\d H=0$, $H\d=0$, so $H$ factors through $M_0\onto\pi_0M$
and over $\pi_1N\rightarrowtail N_1$. 

Hence the proof is finished with the following lemma.
\end{proof}

\begin{Lemma}
For any $B$-modules $M$, $N$ there are isomorphisms
$$
\Ext^n_\a(M,N)\cong\Ext^n_{\pi_0B}(\pi_0M,\pi_0N)
$$
and
$$
(\R^n_\a(\Hom_{\pi_0B}(\pi_0(\_),\pi_1N)))(M)\cong\Ext^n_{\pi_0B}(\pi_0M,\pi_1N).
$$
\end{Lemma}

\begin{proof}
By definition the groups $\Ext^*_\a(M,N)$, respectively
$(\R^n_\a(\Hom_{B_0}(\pi_0(\_),\pi_1N)))(M)$, are cohomology groups of the
complex $[F_\bullet,N]$, resp. $\Hom_{\pi_0B}(\pi_0(F_\bullet),\pi_1N)$,
where $F_\bullet$ is some $\a$-resolution of $M$. We can choose for
$\F_\bullet$ some secondary $\b$-resolution of $M$. Then $\pi_0F_\bullet$
is a free $\pi_0B$-resolution of $\pi_0M$, which makes evident the second
isomorphism. For the first, just note in addition that by \ref{freehom}
$[F_\bullet,N]$ is isomorphic to $\Hom_{B_0}(\pi_0(F_\bullet),\pi_0N)$.
\end{proof}

\section{The track category of spectra}

In this section we introduce the notion of stable maps and stable tracks
between spectra. This yields the track category of spectra. See also
\cite{Baues}*{section 2.5}.

\begin{Definition}\label{sp}
A \emph{spectrum} $X$ is a sequence of maps
$$
X_i\xto r\Omega X_{i+1},\ i\in\Z
$$
in the category $\Topt$ of pointed spaces. This is an $\Omega$-spectrum if
$r$ is a homotopy equivalence for all $i$.

A \emph{stable homotopy class} $f:X\to Y$ between spectra is a sequence of
homotopy classes $f_i\in[X_i,Y_i]$ such that the squares
$$
\xymatrix{
X_i\ar[r]^{f_i}\ar[d]^r&Y_i\ar[d]^r\\
\Omega X_{i+1}\ar[r]^{\Omega f_{i+1}}&\Omega Y_{i+1}
}
$$
commute in $\Topt\ho$. The category $\Sp$ consists of spectra and stable
homotopy classes as morphisms. Its full subcategory $\Omega$-$\Sp$
consisting of $\Omega$-spectra is equivalent to the usual homotopy category
of spectra considered as a Quillen model category.

A \emph{stable map} $f=(f_i,\tilde f_i)_i:X\to Y$ between spectra is a sequence
of diagrams in the track category $\hog{\Topt}$ $(i\in\Z)$
$$
\xymatrix{
X_i\ar[r]^{f_i}\ar[d]_r\drtwocell\omit{^{\tilde f_i\ }}&Y_i\ar[d]^r\\
\Omega X_{i+1}\ar[r]_{\Omega f_{i+1}}&\Omega Y_{i+1}.
}
$$
Obvious composition of such maps yields the category
$$
\hog{\Sp}_0.
$$
It is the underlying category of a track category $\hog{\Sp}$ with tracks
$(H:f\then g)\in\hog{\Sp}_1$ given by sequences
$$
H_i:f_i\then g_i
$$
of tracks in $\Topt$ such that the diagrams
$$
\xymatrix@C=8em@!R=3em{
X_i\ar[r]_{f_i}\ar[d]_r\drtwocell\omit{^{\tilde
f_i\ }}\ruppertwocell^{g_i}{^H_i\ }
&Y_i\ar[d]^r\\
\Omega X_{i+1}\ar[r]^{\Omega f_{i+1}}\rlowertwocell_{\Omega
g_{i+1}}{_\Omega H_{i+1}\hskip4em}
&\Omega Y_{i+1}
}
$$
paste to $\tilde g_i$. This yields a well-defined track category
$\hog{\Sp}$. Moreover
$$
\hog{\Sp}\ho\cong\Sp
$$
is an isomorphism of categories. Let $\hog{X,Y}$ be the groupoid of
morphisms $X\to Y$ in $\hog{\Sp}_0$ and let $\hog{X,Y}_1^0$ be the set of
pairs $(f,H)$ where $f:X\to Y$ is a map and $H:f\then0$ is a track in
$\hog{\Sp}$, i.~e. a stable homotopy class of nullhomotopies for $f$.

For a spectrum $X$ let $\Sigma^kX$ be the \emph{shifted spectrum} with
$(\Sigma^kX)_n=X_{n+k}$ and the commutative diagram
$$
\xymatrix{
(\Sigma^kX)_n\ar@{=}[d]\ar[r]^-r&\Omega(\Sigma^kX)_{n+1}\ar@{=}[d]\\
X_{n+k}\ar[r]^-r&\Omega(X_{n+k+1})
}
$$
defining $r$ for $\Sigma^kX$. A map $f:Y\to\Sigma^kX$ is also called a map
$f$ \emph{of degree $k$} from $Y$ to $X$.
\end{Definition}

\section{The pair algebra $\aB$ and secondary cohomology of spectra as a
$\aB$-module}

The secondary cohomology of a space was introduced in \cite{Baues}*{section
6.3}. We here consider the corresponding notion of secondary cohomology of
a spectrum. 

Let $\F$ be a prime field $\F=\Z/p\Z$ and let $Z$ denote the Eilenberg-Mac
Lane spectrum with 
$$
Z^n=K(\F,n)
$$
chosen as in \cite{Baues}. Here $Z^n$ is a topological $\F$-vector space
and the homotopy equivalence $Z^n\to\Omega Z^{n+1}$ is $\F$-linear. This
shows that for a spectrum $X$ the sets $\hog{X,\Sigma^kZ}_0$ and
$\hog{X,\Sigma^kZ}_1^0$, of stable maps and stable 0-tracks repectively, are $\F$-vector spaces.

We now recall the definition of the pair algebra $\aB=(\d:\aB_1\to\aB_0)$
of secondary cohomology operations from \cite{Baues}. Let $\G=\Z/p^2\Z$ and let
$$
\aB_0=T_\G(E_\A)
$$
be the $\G$-tensor algebra generated by the subset
$$
E_\A=
\begin{cases}
\set{\Sq^1,\Sq^2,...}&\textrm{for $p=2$},\\
\set{\P^1,\P^2,...}\cup\set{\beta,\beta\P^1,\beta\P^2,...}&\textrm{for odd
$p$}
\end{cases}
$$
of the mod $p$ Steenrod algebra $\A$. We define $\aB_1$ by the
pullback diagram of graded abelian groups
\begin{equation}\label{defb1}
\begin{aligned}
\xymatrix{
&\Sigma\A\ar@{ >->}[d]\\
\aB_1\ar[r]\ar[d]_\d\ar@{}[dr]|<{\Bigg{\lrcorner}}&\hog{Z,\Sigma^*Z}_1^0\ar[d]^\d\\
\aB_0\ar[r]^-s&\hog{Z,\Sigma^*Z}_0\ar@{->>}[d]\\
&\A.
}
\end{aligned}
\end{equation}
in which the right hand column is an exact sequence. Here we choose for
$\alpha\in E_\A$ a stable map $s(\alpha):Z\to\Sigma^{|\alpha|}Z$
representing $\alpha$ and we define $s$ to be the $\G$-linear map given on
monomials $a_1\cdots a_n$ in the free monoid Mon$(E_\A)$ generated by
$E_\A$ by the composites
$$
s(a_1\cdots a_n)=s(a_1)\cdots s(a_n).
$$
It is proved in \cite{Baues}*{5.2.3} that $s$ defines a pseudofunctor, that is,
there is a well-defined track
$$
\Gamma:s(a\cdot b)\then s(a)\circ s(b)
$$
for $a,b\in\aB_0$ such that for any $a$, $b$, $c$ pasting of tracks in the
diagram
$$
\includegraphics{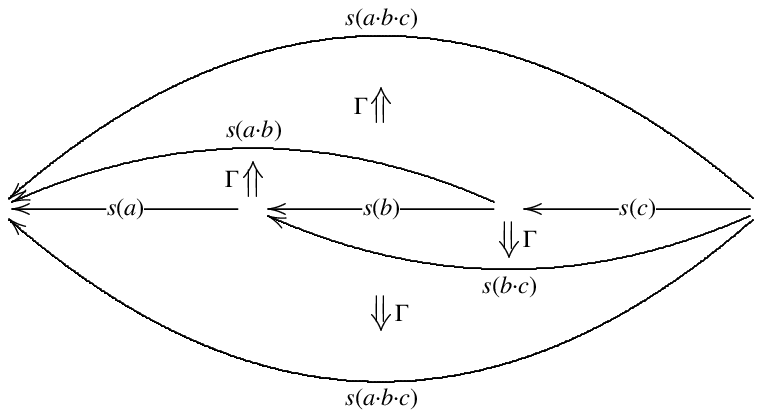}
$$
yields the identity track. Now $\aB_1$ is a $\aB_0$-$\aB_0$-bimodule by
defining
$$
a(b,z)=(a\cdot b,a\bullet z)
$$
with $a\bullet z$ given by pasting $s(a)z$ and $\Gamma$. Similarly
$$
(b,z)a=(b\cdot a,z\bullet a)
$$
where $z\bullet a$ is obtained by pasting $zs(a)$ and $\Gamma$. Then it is
shown in \cite{Baues} that $\aB=(\d:\aB_1\to\aB_0)$ is a well-defined pair
algebra with $\pi_0\aB=\A$ and $\Sigma$-structure $\pi_1\aB=\Sigma\A$.

For a spectrum $X$ let
$$
\H(X)_0=\aB_0\hog{X,\Sigma^*Z}_0
$$
be the free $\aB_0$-module generated by the graded set
$\hog{X,\Sigma^*Z}_0$. We define $\H(X)_1$ by the pullback diagram
$$
\xymatrix{
&\Sigma H^*X\ar@{ >->}[d]\\
\H(X)_1\ar[r]\ar[d]_\d\ar@{}[dr]|<{\Bigg{\lrcorner}}&\hog{X,\Sigma^*Z}_1^0\ar[d]^\d\\
\H(X)_0\ar[r]^-s&\hog{X,\Sigma^*Z}_0\ar@{->>}[d]\\
&H^*X
}
$$
where $s$ is the $\G$-linear map which is the identity on generators and is
defined on words $a_1\cdots a_n\cdot u$ by the composite $s(a_1)\cdots
s(a_n)s(u)$ for $a_i$ as above and $u\in\hog{X,\Sigma^*Z}_0$. Again $s$ is
a pseudofunctor and with actions $\bullet$ defined as above we see that the
graded pair module
$$
\H(X)=\left(\H(X)_1\xto\d\H(X)_0\right)
$$
is a $\aB$-module. We call $\H(X)$ the \emph{secondary cohomology} of the
spectrum $X$. Of course $\H(X)$ has a $\Sigma$-structure in the sense of
\ref{sigma} above.

\begin{Example}
Let $\G^\Sigma$ be the $\aB$-module given by the augmentation $\aB\to\G^
\Sigma$ in \cite{Baues}. Recall that $\G^\Sigma$ is the pair
$$
\G^\Sigma=\left(\F\oplus\Sigma\F\xto\d\G\right)
$$
with $\restr\d\F$ the inclusion nad $\restr\d{\Sigma\F}=0$. Then the
sphere spectrum $S^0$ admits a weak equivalence of $\aB$-modules
$$
\H(S^0)\xto\sim\G^\Sigma.
$$
Compare \cite{Baues}*{12.1.5}.
\end{Example}

\section{The ${\mathrm E}_3$-term of the Adams spectral sequence}

We now are ready to formulate our main result describing the algebraic
equivalent of the E$_3$-term of the Adams spectral sequence. Let $X$ be a
spectrum of finite type and $Y$ a finite dimensional spectrum. Then for
each prime $p$ there is a spectral sequence $\E_*=\E_*(Y,X)$ with
\begin{align*}
\E_*&\Longrightarrow[Y,\Sigma^*X]_p\\
\E_2&=\Ext_\A(H^*X,H^*Y).
\end{align*}

\begin{Theorem}\label{e3der}
The $\E_3$-term $\E_3=\E_3(Y,X)$ of the Adams spectral sequence is given by
the secondary $\Ext$ group defined in \ref{secodef}
$$
\E_3=\Ext_\aB(\H^*X,\H^*Y).
$$
\end{Theorem}

\begin{Corollary}
If $X$ and $Y$ are both the sphere spectrum we get
$$
\E_3(S^0,S^0)=\Ext_\aB(\G^\Sigma,\G^\Sigma).
$$
\end{Corollary}

Since the pair algebra $\aB$ is computed in \cite{Baues} completely we see
that $\E_3(S^0,S^0)$ is algebraically determined. This leads to the
algorithm below computing $\E_3(S^0,S^0)$.

The proof of \ref{e3der} is based on the following result in \cite{Baues}.
Consider the track categories
\begin{align*}
\b&\subset\hog{\Sp}\\
\b'&\subset(\aB\mathrm{-}\Mod)^\Sigma
\end{align*}
where $\hog{\Sp}$ is the track category of spectra in \ref{sp} and
$(\aB\mathrm{-}\Mod)^\Sigma$ is the track category of $\aB$-modules with
$\Sigma$-structure in \ref{sigma} with the pair algebra $\aB$ defined by
\eqref{defb1}. Let $\b$ be the full track subcategory of $\hog{\Sp}$
consisting of finite products of shifted Eilenberg-Mac Lane
spectra $\Sigma^kZ^*$. Moreover let $\b'$ be the full track subcategory of
$(\aB\mathrm{-}\Mod)^\Sigma$ consisting of finitely generated free
$\aB$-modules. As in \cite{BJ5}*{4.3} we obtain for spectra $X$, $Y$ in
\ref{e3der} the track categories
\begin{align*}
\set{Y,X}\b&\subset\hog{\Sp}\\
\b'\set{\H X,\H Y}&\subset(\aB\mathrm{-}\Mod)^\Sigma
\end{align*}
with $\set{Y,X}\b$ obtained by adding to $\b$ the objects $X$, $Y$ and all
morphisms and tracks from $\hog{X,Z}$, $\hog{Y,Z}$ for all objects $Z$ in
$\b$. It is proved in \cite{Baues}*{5.5.6} that the following result holds
which shows that we can apply \cite{BJ5}*{5.1}.

\begin{Theocite}\label{strictif}
There is a strict track equivalence
$$
(\set{Y,X}\b)\op\xto\sim\b'\set{\H X,\H Y}.
$$
\end{Theocite}
\qed

\begin{proof}[Proof of \ref{e3der}]
By the main result 7.3 in \cite{BJ5} we have a description of the
differential $d_{(2)}$ in the Adams spectral sequence by the following
commutative diagram 
$$
\xymatrix{
\Ext^n_{\a\op}(X,Y)^m\ar[r]^{d_{(2)}}\ar[d]^\cong
&\Ext^{n+2}_{\a\op}(X,Y)^{m+1}\ar[d]^\cong\\
\Ext_{\mathscr A}^n(H^*X,H^*Y)^m\ar[r]^{d_{(2)}}
&\Ext_{\mathscr A}^{n+2}(H^*X,H^*Y)^{m+1}
}
$$
where $\a=\b\ho$. On the other hand the differential
$d_{(2)}$ defining the secondary $\Ext$-group $\Ext_\aB(\H X,\H Y)$ is by \ref{sigmacoinc}
given by the commutative diagram
$$
\xymatrix{
\Ext^n_{\a'}(\H X,\H Y)^m\ar[r]\ar@{=}[d]
&\Ext^{n+2}_{\a'}(\H X,\H Y)^{m+1}\ar@{=}[d]\\
\Ext_{\mathscr A}^n(H^*X,H^*Y)^m\ar[r]
&\Ext_{\mathscr A}^{n+2}(H^*X,H^*Y)^{m+1}
}
$$
where $\a'=\b'\ho$. Now \cite{BJ5}*{5.1} shows by \ref{strictif} that the top
rows of these diagrams coincide.
\end{proof}

\section{The multiplication map of $\aB$}

We recall notation $\G=\ZZ/4\ZZ$, $\F=\F_2=\ZZ/2\ZZ$ from \cite{Baues}.
The quotient homomorphism $\G\onto\F$ will be denoted by $\pi$ and the
isomorphism $\F\cong2\G$ by $i$. Moreover we will need the set-theoretic
section $\chi:\F\into\G$ of $\pi$ given by $\chi(0)=0$, $\chi(1)=1$. In
the pair algebra $\aB=\left(\d:\aB_1\to\aB_0\right)$,
recall that $\aB_0$ is the graded free associative $\G$-algebra on the
generators $\Sq^n$ of degree $n$, for $n\ge1$; there is thus a surjective
homomorphism of graded algebras $\pi:\aB_0\onto\A$ onto the mod
2 Steenrod algebra. Its kernel is denoted by $R$, so that we have the short
exact sequence
$$
0\to R\to\aB_0\to\A\to0.
$$
It is well known that $R$ is a graded two-sided ideal generated (as a
two-sided ideal, i.~e. as a $\aB_0$-$\aB_0$-bimodule) by $2\aB_0$ and by
the \emph{Adem elements}
$$
[a,b]:=\Sq^a\Sq^b+\sum_{k=\max(0,a-b+1)}^{\min\left(b-1,\left[\frac
a2\right]\right)}\chi\binom{b-k-1}{a-2k}\Sq^{a+b-k}\Sq^k,
$$
for $0<a<2b$. As shown in \cite{Baues}, one can generate $R$ as a right
$\aB_0$-module by $2\in R^0$ and the \emph{admissible relations}
$\Sq^{a_k}\Sq^{a_{k-1}}\cdots\Sq^{a_1}[a_0,b]\in R^{a_k+...+a_0+b}$, with
$a_{j+1}\ge2a_j$ for all $k>j\ge0$ and $a_0<2b$.

As for the rest of the structure of $\aB$, as an abelian group,
$\aB_1$ is $R\oplus\Sigma\A$, that is,
$$
\aB_1^n=R^n\oplus\A^{n-1},
$$
and $\d$ is the projection. Moreover the $\aB_0$-bimodule
structure of $\aB_1$ is given by
$$
(r,a)b=(rb,a\pi(b))
$$
and
$$
b(r,a)=(br,A(\pi(b),r)+\pi(b)a),
$$
where
$$
A:\A\ox R\to\A
$$
is the \emph{multiplication map} of degree -1 described in \cite{Baues}.
Algebraic properties characterizing the multiplication map $A$ are
achieved in \cite{Baues}*{theorem 16.3.3}. In \cite{Baues}*{section 16.6}
an algorithm is obtained which computes the multiplication map $A$.

Particular important elements of $\aB_1$ get special notation; e.~g. we
have $[2]:=(2,0)\in R^0\oplus0=\aB_1^0$ and $\Sigma1:=(0,1)\in
R^1\oplus\A^0=\aB_1^1$. This pair algebra has an augmentation
$\epsilon:\aB\to\G^\Sigma$, where
$\G^\Sigma=\left((i,0):\F\oplus\Sigma\F\to\G\right)$ is the graded
$\aB$-module equal to $i:\F\cong2\G\subset\G$ in degree 0, to $\F\to0$ in
degree 1 and zero in all other degrees. Components of $\epsilon$ are the
augmentation $\epsilon_0:\aB_0\to\G$ and the homomorphism
$\epsilon_1:R\oplus\Sigma\A\to2\G\oplus\Sigma\F$ given by
$(r,a)\mapsto(\epsilon_0(r),\epsilon(a))$.

\section{The algorithm for the computation of $d_{(2)}$ on
$\Ext_\A(\F,\F)$ in terms of the multiplication maps}

Suppose now given some projective resolution of the left $\mathscr
A$-module $\F$. For definiteness, we will work with the minimal
resolution
\begin{equation}\label{minireso}
\F\ot\A\brk{g_0^0}\ot\A\brk{g_1^{2^n}\mid
n\ge0}\ot\A\brk{g_2^{2^i+2^j}\mid
|i-j|\ne1}\ot...,
\end{equation}
where $g_m^d$, $d\ge m$, is a generator of the $m$-th resolving module in degree $d$.
Sometimes there are more than one generators with the same $m$ and $d$, in which
case the further ones will be denoted by $'g_m^d$, $''g_m^d, \cdots$.

These generators and values of the differential on them can be computed
effectively; for example, $d(g_1^{2^n})=\Sq^{2^n}g_0^0$ and
$d(g_m^m)=\Sq^1g_{m-1}^{m-1}$; moreover e.~g. an
algorithm from \cite{Bruner} gives
$$
\alignbox{
d(g_2^4)&=\Sq^3g_1^1+\Sq^2g_1^2\\
d(g_2^5)&=\Sq^4g_1^1+\Sq^2\Sq^1g_1^2+\Sq^1g_1^4\\
d(g_2^8)&=\Sq^6g_1^2+(\Sq^4+\Sq^3\Sq^1)g_1^4\\
d(g_2^9)&=\Sq^8g_1^1+(\Sq^5+\Sq^4\Sq^1)g_1^4+\Sq^1g_1^8\\
d(g_2^{10})&=(\Sq^8+\Sq^5\Sq^2\Sq^1)g_1^2+(\Sq^5\Sq^1+\Sq^4\Sq^2)g_1^4+\Sq^2g_1^8\\
d(g_2^{16})&=(\Sq^{12}+\Sq^9\Sq^2\Sq^1+\Sq^8\Sq^3\Sq^1)g_1^4+(\Sq^8+\Sq^7\Sq^1+\Sq^6\Sq^2)g_1^8\\
\cdots,\\
d(g_3^6)&=\Sq^4g_2^2+\Sq^2g_2^4+\Sq^1g_2^5\\
d(g_3^{10})&=\Sq^8g_2^2+(\Sq^5+\Sq^4\Sq^1)g_2^5+\Sq^1g_2^9\\
d(g_3^{11})&=(\Sq^7+\Sq^4\Sq^2\Sq^1)g_2^4+\Sq^6g_2^5+\Sq^2\Sq^1g_2^8\\
d(g_3^{12})&=\Sq^8g_2^4+(\Sq^6\Sq^1+\Sq^5\Sq^2)g_2^5+(\Sq^4+\Sq^3\Sq^1)g_2^8+\Sq^3g_2^9+\Sq^2g_2^{10}\\
\cdots,
\\
d(g_4^{11})&=\Sq^8g_3^3+(\Sq^5+\Sq^4\Sq^1)g_3^6+\Sq^1g_3^{10}\\
d(g_4^{13})&=\Sq^8\Sq^2g_3^3+(\Sq^7+\Sq^4\Sq^2\Sq^1)g_3^6+\Sq^2\Sq^1g_3^{10}+\Sq^2g_3^{11}\\
\cdots,
\\
d(g_5^{14})&=\Sq^{10}g_4^4+\Sq^2\Sq^1g_4^{11}\\
d(g_5^{16})&=\Sq^{12}g_4^4+\Sq^4\Sq^1g_4^{11}+\Sq^3g_4^{13}\\
\cdots,
\\
d(g_6^{16})&=\Sq^{11}g_5^5+\Sq^2g_5^{14}\\
\cdots,
}
$$
etc.

By understanding the above formul\ae\ \emph{literally} (i.~e. by applying
$\chi$ degreewise to them), each such resolution gives rise to a sequence
of $\aB$-module homomorphisms
\begin{equation}\label{precomplex}
\G^\Sigma\ot\aB\brk{g_0^0}\ot\aB\brk{g_1^{2^n}\mid
n\ge0}\ot\aB\brk{g_2^{2^i+2^j}\mid
|i-j|\ne1}\ot...,
\end{equation}
which is far from being exact --- in fact even the composites of consecutive
maps are not zero. In more detail, one has commutative diagrams
$$
\xymatrix{
2\G\ar[d]
&R^0g_0^0\ar[l]_{\epsilon_0}\ar[d]
&0\ar[l]\ar[d]
&...\ar[l]\\
\G
&\aB_0^0g_0^0\ar[l]_{\epsilon_0}
&0\ar[l]
&...\ar[l]
}
$$
in degree 0,
$$
\xymatrix{
\F\ar[d]
&R^1g_0^0\oplus\A^0g_0^0\ar[l]_-{(0,\epsilon)}\ar[d]
&R^0g_1^1\ar[l]_-{\binom d0}\ar[d]
&0\ar[l]\ar[d]
&...\ar[l]\\
0
&\aB_0^1g_0^0\ar[l]
&\aB_0^0g_1^1\ar[l]_d
&0\ar[l]
&...\ar[l]
}
$$
in degree 1,
$$
\xymatrix{
0\ar[d]
&R^2g_0^0\oplus\A^1g_0^0\ar[l]\ar[d]
&\left(R^1g_1^1\oplus R^0g_1^2\right)\oplus{\mathscr
A}^0g_1^1\ar[l]_-{\smat{d&0\\0&d}}\ar[d]
&R^0g_2^2\ar[l]_-{\binom d0}\ar[d]
&0\ar[l]\ar[d]
&...\ar[l]\\
0
&\aB_0^2g_0^0\ar[l]
&\aB_0^1g_1^1\oplus\aB_0^0g_1^2\ar[l]_-d
&\aB_0^0g_2^2\ar[l]_-d
&0\ar[l]
&...\ar[l]
}
$$
in degree 2, ...
$$
\xymatrix{
0\ar[d]
&R^ng_0^0\oplus\A^{n-1}g_0^0\ar[l]\ar[d]
&\bigoplus_{2^i\le n}R^{n-2^i}g_1^{2^i}
\oplus\bigoplus_{2^i\le n-1}\A^{n-1-2^i}g_1^{2^i}\ar[l]_-{\smat{d&0\\0&d}}\ar[d]
&...\ar[l]\\
0
&\aB_0^ng_0^0\ar[l]
&\bigoplus_{2^i\le n}\aB_0^{n-2^i}g_1^{2^i}\ar[l]_-d
&...\ar[l]
}
$$
in degree $n$, etc.

Our task is then to complete these diagrams into an exact secondary complex
via certain (degree preserving) maps 
$$
\delta_m=\binom{\delta^R_m}{\delta^\A_m}:\aB_0\brk{g_{m+2}^n\mid n}\to(R\oplus\Sigma\A)\brk{g_m^n\mid
n}.
$$

Now for these maps to form a secondary complex, according to
\ref{secs}.1 one must have $\d\delta=d_0d_0$,
$\delta\d=d_1d_1$, and $d_1\delta=\delta d_0$. One sees easily that
these equations together with the requirement that $\delta$ be left
$\aB_0$-module homomorphism are equivalent to 
\begin{align}
\delta^R&=dd,\\
\label{deltaeqs}\delta^\A(bg)&=\pi(b)\delta^\A(g)+A(\pi(b),dd(g)),\\
d\delta^\A&=\delta^\A d,
\end{align}
for $b\in\aB_0$, $g$ one of the $g_m^n$, and $A(a,rg):=A(a,r)g$ for
$a\in\A$, $r\in R$. Hence $\delta$ is completely determined by the elements
$$
\delta^\A_m(g_{m+2}^n)\in\bigoplus_k\A^{n-k-1}\brk{g_m^k}
$$
which, to form a secondary complex, are only required to satisfy
$$
d\delta_m^\A(g_{m+2}^n)=\delta_{m-1}^\A d(g_{m+2}^n),
$$
where on the right $\delta_{m-1}^\A$ is extended to
$\aB_0\brk{g_{m+1}^*}$ via \ref{deltaeqs}. Then furthermore secondary
exactness must hold, which by \ref{secs} means that the (ordinary) complex
$$
\ot
\aB_0\brk{g_{m-1}^*}\oplus(R\oplus\Sigma\A)\brk{g_{m-2}^*}
\ot
\aB_0\brk{g_m^*}\oplus(R\oplus\Sigma\A)\brk{g_{m-1}^*}
\ot
\aB_0\brk{g_{m+1}^*}\oplus(R\oplus\Sigma\A)\brk{g_m^*}
\ot
$$
with differentials
$$
\smat{d_{m+1}&i_{m+1}&0\\d_md_{m+1}&d_m&0\\\delta_m^\A&0&d_m}:
\aB_0\brk{g_{m+2}^*}\oplus R\brk{g_{m+1}^*}\oplus\Sigma\A\brk{g_{m+1}^*}
\to
\aB_0\brk{g_{m+1}^*}\oplus R\brk{g_m^*}\oplus\Sigma\A\brk{g_m^*}
$$
is exact. Then straightforward checking shows that one can eliminate $R$
from this complex altogether, so that its exactness is equivalent to the exactness
of a smaller complex
$$
\ot
\aB_0\brk{g_{m-1}^*}\oplus\Sigma\A\brk{g_{m-2}^*}
\ot
\aB_0\brk{g_m^*}\oplus\Sigma\A\brk{g_{m-1}^*}
\ot
\aB_0\brk{g_{m+1}^*}\oplus\Sigma\A\brk{g_m^*}
\ot
$$
with differentials
$$
\smat{d_{m+1}&0\\\delta_m^\A&d_m}:
\aB_0\brk{g_{m+2}^*}\oplus\Sigma\A\brk{g_{m+1}^*}
\to
\aB_0\brk{g_{m+1}^*}\oplus\Sigma\A\brk{g_m^*}.
$$
Note also that by \ref{deltaeqs} $\delta^\A$ factors through
$\pi$ to give
$$
\bar\delta_m:\A\brk{g^*_{m+2}}\to\Sigma\A\brk{g^*_m}.
$$
It follows that secondary exactness of the resulting complex is equivalent
to exactness of the \emph{mapping cone} of this $\bar\delta$, i.~e. to the
requirement that $\bar\delta$ is a quasiisomorphism. On the other hand, the
complex $(\A\brk{g^*_*},d_*)$ is acyclic by construction, so any
of its self-maps is a quasiisomorphism. We thus obtain

\begin{Theorem}
Completions of the diagram \ref{precomplex}
to an exact secondary complex are in one-to-one correspondence with
maps $\delta_m:\A\brk{g_{m+2}^*}\to\Sigma\A\brk{g^*_m}$
satisfying
\begin{equation}\label{maineq}
d\delta g=\delta dg,
\end{equation}
with $\delta(ag)$ for $a\in\A$ defined by
$$
\delta(ag)=a\delta(g)+A(a,ddg)
$$
where $A(a,rg)$ for $r\in R$ is interpreted as $A(a,r)g$.
\end{Theorem}
\qed

We can use this to construct the secondary resolution inductively.
Just start by introducing values of $\delta$ on the generators as
expressions with indeterminate coefficients; the equation \eqref{maineq} will
impose linear conditions on these coefficients. These are then solved
degree by degree. For example, in degree 2 one may have
$$
\delta(g_2^2)=\eta_2^2(\Sq^1)\Sq^1g_0^0
$$
for some $\eta_2^2(\Sq^1)\in\F$. Similarly in degree 3 one may have
$$
\delta(g_3^3)=\eta_3^3(\Sq^1)\Sq^1g_1^1+\eta_3^3(1)g_1^2.
$$
Then one will get
$$
d\delta(g_3^3)=\eta_3^3(\Sq^1)\Sq^1d(g_1^1)+\eta_3^3(1)d(g_1^2)=\eta_3^3(\Sq^1)\Sq^1\Sq^1g_0^0+\eta_3^3(1)\Sq^2g_0^0=\eta_3^3(1)\Sq^2g_0^0
$$
and
\begin{multline*}
\delta d(g_3^3)=\delta(\Sq^1g_2^2)\\
=\Sq^1\delta(g_2^2)+A(\Sq^1,dd(g_2^2))=\eta_2^2(\Sq^1)\Sq^1\Sq^1g_0^0+A(\Sq^1,d(\Sq^1g_1^1))=A(\Sq^1,\Sq^1\Sq^1g_0^0)=0;
\end{multline*}
thus \eqref{maineq} forces $\eta^3_3(1)=0$.

Similarly one puts $\delta(g_m^d)=\sum_{m-2\le d'\le d-1}\sum_a\eta_m^d(a)ag_{m-2}^{d'}$,
with $a$ running over a basis in $\A^{d-1-d'}$, and then substituting this in
\eqref{maineq} gives linear equations on the numbers $\eta_m^d(a)$. Solving
these equations and choosing the remaining $\eta$'s arbitrarily then gives
values of the differential $\delta$ in the secondary resolution.

Then finally to obtain the secondary differential
$$
d_{(2)}:\Ext^n_\A(\F,\F)^m\to\Ext^{n+2}_\A(\F,\F)^{m+1}
$$
from this $\delta$, one just applies the functor $\Hom_\A(\_,\F)$ to the
initial minimal resolution and calculates the map induced by $\delta$ on
cohomology of the resulting cochain complex, i.~e. on $\Ext^*_\A(\F,\F)$.
In fact since \eqref{minireso} is a minimal resolution, the value of
$\Hom_\A(\_,\F)$ on it coincides with its own cohomology and is the
$\F$-vector space of those linear maps $\A\brk{g_*^*}\to\F$ which vanish
on all elements of the form $ag_*^*$ with $a$ of positive degree.

Let us then identify $\Ext^*_\A(\F,\F)$ with this space and choose a basis
in it consisting of elements $\hat g_m^d$ defined as the maps sending the
generator $g_m^d$ to 1 and all other generators to 0. One then has
$$
(d_{(2)}(\hat g_m^d))(g_{m'}^{d'})=\hat g_m^d\delta(g_{m'}^{d'}).
$$ 
The right hand side is nonzero precisely when $g_m^d$ appears in
$\delta(g_{m'}^{d'})$ with coefficient 1, i.~e. one has
$$
d_{(2)}(\hat g_m^d)=\sum_{\textrm{$g_m^d$ appears in
$\delta(g_{m+2}^{d+1})$}}\hat g_{m+2}^{d+1}.
$$

For example, looking at the table of values of $\delta$ below we see that
the first instance of a $g_m^d$ appearing with coefficient 1 in a value of
$\delta$ on a generator is
$$
\delta(g_3^{17})=g_1^{16}+ \Sq^{12} g_1^4+\Sq^{10}\Sq^4
g_1^2+(\Sq^9\Sq^4\Sq^2+\Sq^{10}\Sq^5+\Sq^{11}\Sq^4)g_1^1.
$$
This means
$$
d_{(2)}(\hat g_1^{16})=\hat g_3^{17}
$$
and moreover $d_{(2)}(\hat g_m^d)=0$ for all $g_m^d$ with $d<17$ (one can
check all cases for each given $d$ since the number of generators $g_m^d$
for each given $d$ is finite).

Treating similarly the rest of the table below we find that the only
nonzero values of $d_{(2)}$ on generators of degree $<36$ are as follows:
$$
\begin{array}{rl}
d_{(2)}(\hat g_1^{16})&=\hat g_3^{17}\\
d_{(2)}(\hat g_4^{21})&=\hat g_6^{22}\\
d_{(2)}(\hat g_4^{22})&=\hat g_6^{23}\\
d_{(2)}(\hat g_5^{23})&=\hat g_7^{24}\\
d_{(2)}(\hat g_7^{30})&=\hat g_9^{31}\\
d_{(2)}(\hat g_8^{31})&=\hat g_{10}^{32}\\
d_{(2)}(\hat g_1^{32})&=\hat g_3^{33}\\
d_{(2)}(\hat g_2^{33})&=\hat g_4^{34}\\
d_{(2)}(\hat g_7^{33})&=\hat g_9^{34}\\
d_{(2)}(\hat g_8^{33})&=\hat g_{10}^{34}\\
d_{(2)}({}'\hat g_3^{34})&=\hat g_5^{35}\\
d_{(2)}(\hat g_8^{34})&=\hat g_{10}^{35}.\\
\end{array}
$$
Presently a computer calculation continues reaching degree 39 and showing
that up to that degree there are the following further nonzero
differentials:
$$
\begin{array}{rl}
d_{(2)}(\hat g_7^{36})&=\hat g_9^{37}\\
d_{(2)}(\hat g_8^{37})&=\hat g_{10}^{38}.
\end{array}
$$
These data can be summarized in the following picture, thus confirming
calculations presented in Ravenel's book \cite{Ravenel}.

\ 

\newcount\s
\def\latticebody{
\s=\latticeA
\advance\s by\latticeB
\ifnum\s<40\drop{.}\else\fi}
\ \hfill
\xy
*\xybox{
0;<.3cm,0cm>:
,0,{\xylattice0{40}0{15}}
,(0,0)*{\bullet}
,(0,1)*{\bullet}
,(0,2)*{\bullet}
,(0,3)*{\bullet}
,(0,4)*{\bullet}
,(0,5)*{\bullet}
,(0,6)*{\bullet}
,(0,7)*{\bullet}
,(0,8)*{\bullet}
,(0,9)*{\bullet}
,(0,10)*{\bullet}
,(0,11)*{\bullet}
,(0,12)*{\bullet}
,(0,13)*{\bullet}
,(0,14)*{\bullet}
,(0,15)*{\bullet}
,(1,1)*{\bullet}
,(3,1)*{\bullet}
,(7,1)*{\bullet}
,(15,1)*{\circ}="1.16"
,(31,1)*{\circ}="1.32"
,(2,2)*{\bullet}
,(3,2)*{\bullet}
,(6,2)*{\bullet}
,(7,2)*{\bullet}
,(8,2)*{\bullet}
,(14,2)*{\bullet}
,(15,2)*{\bullet}
,(16,2)*{\bullet}
,(18,2)*{\bullet}
,(30,2)*{\bullet}
,(31,2)*{\circ}="2.33"
,(32,2)*{\bullet}
,(34,2)*{\bullet}
,(3,3)*{\bullet}
,(7,3)*{\bullet}
,(8,3)*{\bullet}
,(9,3)*{\bullet}
,(14,3)*{\circ};"1.16"**\dir{-}
,(15,3)*{\bullet}
,(17,3)*{\bullet}
,(18,3)*{\bullet}
,(19,3)*{\bullet}
,(21,3)*{\bullet}
,(30,3)*{\circ};"1.32"**\dir{-}
,(31,3)*{\circ}="3.34"
,(33,3)*{\bullet}
,(34,3)*{\bullet}
,(7,4)*{\bullet}
,(9,4)*{\bullet}
,(14,4)*{\bullet}
,(15,4)*{\bullet}
,(17,4)*{\circ}="4.21"
,(18.15,3.85)*{\bullet}
,(17.85,4.15)*{\circ}="4.22"
,(20,4)*{\bullet}
,(22,4)*{\bullet}
,(23,4)*{\bullet}
,(30,4)*{\circ};"2.33"**\dir{-}
,(31,4)*{\bullet}
,(32,4)*{\bullet}
,(33,4)*{\bullet}
,(34,4)*{\bullet}
,(9,5)*{\bullet}
,(11,5)*{\bullet}
,(14,5)*{\bullet}
,(14.85,5.15)*{\bullet}
,(15.15,4.85)*{\bullet}
,(17,5)*{\bullet}
,(18,5)*{\circ}="5.23"
,(20,5)*{\bullet}
,(21,5)*{\bullet}
,(23,5)*{\bullet}
,(24,5)*{\bullet}
,(30,5)*{\circ};"3.34"**\dir{-}
,(30.85,5.15)*{\bullet}
,(31.15,4.85)*{\bullet}
,(33,5)*{\bullet}
,(10,6)*{\bullet}
,(11,6)*{\bullet}
,(14,6)*{\bullet}
,(15,6)*{\bullet}
,(16,6)*{\circ};"4.21"**\dir{-}
,(17,6)*{\circ};"4.22"**\dir{-}
,(20,6)*{\bullet}
,(23,6)*{\bullet}
,(26,6)*{\bullet}
,(30,6)*{\bullet}
,(31,6)*{\bullet}
,(32,6)*{\bullet}
,(11,7)*{\bullet}
,(15,7)*{\bullet}
,(16,7)*{\bullet}
,(17,7)*{\circ};"5.23"**\dir{-}
,(23,7)*{\circ}="7.30"
,(26,7)*{\circ}="7.33"
,(29,7)*{\circ}="7.36"
,(30,7)*{\bullet}
,(31,7)*{\bullet}
,(32,7)*{\bullet}
,(15,8)*{\bullet}
,(17,8)*{\bullet}
,(22,8)*{\bullet}
,(23,8)*{\circ}="8.31"
,(25,8)*{\circ}="8.33"
,(26,8)*{\circ}="8.34"
,(28,8)*{\bullet}
,(29,8)*{\circ}="8.37"
,(30,8)*{\bullet}
,(30.85,8.15)*{\bullet}
,(31.15,7.85)*{\bullet}
,(17,9)*{\bullet}
,(19,9)*{\bullet}
,(22,9)*{\circ};"7.30"**\dir{-}
,(22.85,9.15)*{\bullet}
,(23.15,8.85)*{\bullet}
,(25,9)*{\circ};"7.33"**\dir{-}
,(26,9)*{\circ}="9.35"
,(28,9)*{\circ};"7.36"**\dir{-}
,(29,9)*{\bullet}
,(30,9)*{\bullet}
,(18,10)*{\bullet}
,(19,10)*{\bullet}
,(22,10)*{\circ};"8.31"**\dir{-}
,(23,10)*{\bullet}
,(24,10)*{\circ};"8.33"**\dir{-}
,(25,10)*{\circ};"8.34"**\dir{-}
,(28,10)*{\circ};"8.37"**\dir{-}
,(19,11)*{\bullet}
,(23,11)*{\bullet}
,(24,11)*{\bullet}
,(25,11)*{\circ};"9.35"**\dir{-}
,(23,12)*{\bullet}
,(25,12)*{\bullet}
,(25,13)*{\bullet}
,(0,17)*{\ }
,(42,0)*{\ }
}="O"
,{"O"+L \ar "O"+R*+!LD{d-m}}
,{"O"+D \ar "O"+U*+!RD{m}}
\endxy
\hfill\ 

\section{The table of values of the differential $\delta$ in the secondary
resolution for $\G^\Sigma$}

The following table presents results of computer calculations of the
differential $\delta$. Note that it does not have invariant meaning since
it depends on the choices involved in determination of the multiplication
map $A$, of the resolution and of those indeterminate coefficients
$\eta_m^d(a)$ which remain undetermined after the conditions
\eqref{maineq} are satisfied. The resulting secondary differential
$d_{(2)}$ however does not depend on these choices and is canonically
determined.

$$
\begin{array}{rl}
\delta(g_2^2) &= 0\\
\ \\
\delta(g_3^3) &= 0\\
\ \\
\delta(g_2^4) &= 0\\
\delta(g_4^4) &= 0\\
\ \\
\delta(g_2^5) &= 0\\
\delta(g_5^5) &= 0\\
\ \\
\delta(g_3^6) &= \Sq^4 g_1^1
\end{array}
$$

$$
\begin{array}{rl}
\delta(g_6^6) &= 0\\
\ \\
\delta(g_7^7) &= 0\\
\ \\
\delta(g_2^8) &= 0\\
\delta(g_8^8) &= 0\\
\ \\
\delta(g_2^9) &= 0\\
\delta(g_9^9) &= 0\\
\ \\
\delta(g_2^{10}) &= 0\\
\delta(g_3^{10}) &= (\Sq^4\Sq^2\Sq^1 + \Sq^7) g_1^2\\
 &+ \Sq^8 g_1^1\\
\delta(g_{10}^{10}) &= 0\\
\ \\
\delta(g_3^{11}) &= (\Sq^7\Sq^1 + \Sq^8) g_1^2\\
 &+ \Sq^6\Sq^3 g_1^1\\
\delta(g_4^{11}) &= \Sq^5 g_2^5\\
 &+\Sq^4\Sq^2 g_2^4\\
\delta(g_{11}^{11}) &= 0\\
\ \\
\delta(g_3^{12}) &= \Sq^7\Sq^3 g_1^1\\
\delta(g_{12}^{12}) &= 0\\
\ \\
\delta(g_4^{13}) &= \Sq^4 g_2^8\\
 &+ (\Sq^7 + \Sq^5\Sq^2) g_2^5\\
 &+ (\Sq^8 + \Sq^6\Sq^2) g_2^4\\
 &+ (\Sq^7\Sq^3 + \Sq^8\Sq^2 + \Sq^{10}) g_2^2\\
\delta(g_{13}^{13}) &= 0\\
\ \\
\delta(g_5^{14}) &= \Sq^4\Sq^2\Sq^1 g_3^6\\
 &+ (\Sq^7\Sq^3 + \Sq^8\Sq^2) g_3^3\\
\delta(g_{14}^{14}) &= 0\\
\ \\
\delta(g_2^{16}) &= 0\\
\delta(g_5^{16}) &=\Sq^3 g_3^{12}\\
 &+ \Sq^4 g_3^{11}\\
 &+ \Sq^5 g_3^{10}\\
 &+ \Sq^{10}\Sq^2 g_3^3\\ 
\delta(g_6^{16}) &= 0\\
\ \\
\delta(g_2^{17}) &= 0\\
\delta(g_3^{17}) &= g_1^{16}\\
 &+ \Sq^{12} g_1^4\\
 &+\Sq^{10}\Sq^4 g_1^2\\
 &+ (\Sq^9\Sq^4\Sq^2 + \Sq^{10}\Sq^5 + \Sq^{11}\Sq^4) g_1^1\\ 
\delta(g_6^{17}) &= (\Sq^5 + \Sq^4\Sq^1) g_4^{11}\\
 &+ (\Sq^{12} +\Sq^{10}\Sq^2) g_4^4\\
\ \\
\delta(g_2^{18}) &= 0\\
\delta(g_3^{18}) &=(\Sq^{11}\Sq^4 + \Sq^8\Sq^4\Sq^2\Sq^1) g_1^2\\
 &+(\Sq^{10}\Sq^4\Sq^2 + \Sq^{11}\Sq^5 + \Sq^{12}\Sq^4 +\Sq^{14}\Sq^2 + \Sq^{16}) g_1^1
\end{array}
$$

$$
\begin{array}{rl}
\delta(g_4^{18}) &=(\Sq^6\Sq^1 + \Sq^7) g_2^{10}\\
 &+(\Sq^6\Sq^3 + \Sq^7\Sq^2 + \Sq^9) g_2^8\\
 &+\Sq^8\Sq^4g_2^5\\
 &+(\Sq^{10}\Sq^2\Sq^1 + \Sq^{13} + \Sq^{11}\Sq^2 +\Sq^{12}\Sq^1)g_2^4\\
 &+(\Sq^9\Sq^4\Sq^2 + \Sq^{15} + \Sq^{12}\Sq^3 +\Sq^{10}\Sq^5)g_2^2\\
\delta(g_7^{18}) &= \Sq^2\Sq^1 g_5^{14}\\
\ \\
\delta(g_4^{19}) &=\Sq^9 g_2^9\\
 &+(\Sq^{10} + \Sq^8\Sq^2) g_2^8\\
 &+ \Sq^{11}\Sq^2g_2^5\\
 &+ ( \Sq^{11}\Sq^2\Sq^1 + \Sq^{13}\Sq^1 + \Sq^8\Sq^4\Sq^2 + \Sq^{10}\Sq^3\Sq^1) g_2^4\\
 &+(\Sq^{14}\Sq^2 +\Sq^{10}\Sq^4\Sq^2 + \Sq^{12}\Sq^4) g_2^2\\
\delta(g_5^{19}) &=\Sq^1 g_3^{17}\\
 &+ \Sq^4\Sq^2 g_3^{12}\\
 &+ \Sq^4\Sq^2\Sq^1 g_3^{11}\\ 
 &+ (\Sq^6\Sq^2 + \Sq^8) g_3^{10}\\
 &+ (\Sq^8\Sq^4 +\Sq^{11}\Sq^1) g_3^6\\
 &+ (\Sq^{13}\Sq^2 + \Sq^{10}\Sq^5 + \Sq^{15} + \Sq^{11}\Sq^4) g_3^3\\
\ \\
\delta(g_2^{20}) &= 0\\
\delta(g_3^{20}) &= (\Sq^{15} + \Sq^9\Sq^4\Sq^2) g_1^4\\
 &+ (\Sq^{12}\Sq^5 + \Sq^{13}\Sq^4 + \Sq^{16}\Sq^1) g_1^2\\
 &+ (\Sq^{11}\Sq^5\Sq^2 + \Sq^{15}\Sq^3 + \Sq^{18} + \Sq^{12}\Sq^6)g_1^1\\
\delta(g_5^{20}) &= \Sq^4\Sq^2\Sq^1 g_3^{12}\\
 &+ (\Sq^7\Sq^1 + \Sq^8)g_3^{11}\\
 &+ (\Sq^{10}\Sq^3 + \Sq^8\Sq^4\Sq^1 + \Sq^{13} + \Sq^{11}\Sq^2)g_3^6\\
 &+ (\Sq^{13}\Sq^3 + \Sq^{10}\Sq^4\Sq^2 + \Sq^{11}\Sq^5 + \Sq^{12}\Sq^4)g_3^3\\
\delta({}'g_5^{20}) &= \Sq^5\Sq^2 g_3^{12}\\
 &+ \Sq^7\Sq^2 g_3^{10}\\
 &+ (\Sq^{12}\Sq^1 + \Sq^{10}\Sq^3 + \Sq^8\Sq^4\Sq^1 + \Sq^{10}\Sq^2\Sq^1
 + \Sq^{11}\Sq^2) g_3^6\\
 &+ (\Sq^{14}\Sq^2 + \Sq^{13}\Sq^3 + \Sq^{11}\Sq^5 + \Sq^{16} + \Sq^{12}\Sq^4)g_3^3\\
\delta(g_6^{20}) &= (\Sq^6\Sq^2 + \Sq^8) g_4^{11}\\
 &+ (\Sq^{13}\Sq^2 + \Sq^{15} + \Sq^{11}\Sq^4) g_4^4\\
\ \\
\delta(g_3^{21}) &= (\Sq^{15}\Sq^2\Sq^1 + \Sq^{17}\Sq^1 +
 \Sq^{12}\Sq^6) g_1^2\\
 &+ ( \Sq^{13}\Sq^4\Sq^2 + \Sq^{15}\Sq^4 + \Sq^{16}\Sq^3 + \Sq^{17}\Sq^2 + \Sq^{19}) g_1^1\\
\delta(g_4^{21}) &=\Sq^3 g_2^{17}\\
 &+ (\Sq^{10} +\Sq^9\Sq^1) g_2^{10}\\
 &+ (\Sq^9\Sq^3 + \Sq^{11}\Sq^1) g_2^8\\
 &+ (\Sq^{15} +\Sq^{13}\Sq^2 + \Sq^{10}\Sq^5) g_2^5\\
 &+ ( \Sq^{13}\Sq^2\Sq^1 + \Sq^{12}\Sq^3\Sq^1 + \Sq^{12}\Sq^4 + \Sq^9\Sq^4\Sq^2\Sq^1 + \Sq^{10}\Sq^4\Sq^2) g_2^4\\
 &+ (\Sq^{16}\Sq^2 + \Sq^{12}\Sq^6 + \Sq^{15}\Sq^3) g_2^2\\
\delta(g_6^{21}) &=  (\Sq^7 + \Sq^6\Sq^1) g_4^{13}\\
 &+ (\Sq^9 + \Sq^8\Sq^1) g_4^{11}\\
 &+ \Sq^{11}\Sq^5 g_4^4\\
\ \\
\delta(g_3^{22}) &=\Sq^{17} g_1^4\\
 &+ (\Sq^{16}\Sq^2\Sq^1 + \Sq^{13}\Sq^6 +\Sq^{12}\Sq^4\Sq^2\Sq^1 + \Sq^{12}\Sq^6\Sq^1) g_1^2\\
 &+ ( \Sq^{13}\Sq^5\Sq^2 + \Sq^{17}\Sq^3 + \Sq^{18}\Sq^2 +\Sq^{14}\Sq^4\Sq^2) g_1^1
\end{array}
$$

$$
\begin{array}{rl}
\delta(g_4^{22}) &=\Sq^4 g_2^{17}\\
 &+ \Sq^{11} g_2^{10}\\
 &+ (\Sq^{12} + \Sq^9\Sq^3) g_2^9\\
 &+ (\Sq^9\Sq^4 + \Sq^{13} + \Sq^8\Sq^4\Sq^1)g_2^8\\
 &+ \Sq^{12}\Sq^4g_2^5\\
 &+ \Sq^{15}\Sq^2 g_2^4\\
 &+ (\Sq^{13}\Sq^4\Sq^2 + \Sq^{19} + \Sq^{13}\Sq^6 + \Sq^{14}\Sq^5)g_2^2\\
\delta({}'g_4^{22}) &= \Sq^2\Sq^1 g_2^{18}\\
 &+ (\Sq^8\Sq^4 + \Sq^{12}) g_2^9\\
 &+ (\Sq^9\Sq^4 + \Sq^{13} + \Sq^{12}\Sq^1) g_2^8\\
 &+ (\Sq^{16} + \Sq^{13}\Sq^3) g_2^5\\
 &+ (\Sq^{15}\Sq^2 + \Sq^{16}\Sq^1 + \Sq^{13}\Sq^4 + \Sq^{11}\Sq^4\Sq^2)g_2^4\\
 &+ (\Sq^{14}\Sq^5 + \Sq^{19} + \Sq^{17}\Sq^2) g_2^2\\
\delta(g_5^{22}) &= (\Sq^7\Sq^2 + \Sq^6\Sq^2\Sq^1 + \Sq^6\Sq^3)g_3^{12}\\
 &+ \Sq^{10} g_3^{11} + (\Sq^9\Sq^2 + \Sq^8\Sq^3 + \Sq^{11}) g_3^{10}\\
 &+(\Sq^{14}\Sq^1 + \Sq^{11}\Sq^3\Sq^1 + \Sq^{12}\Sq^3 +\Sq^{13}\Sq^2)g_3^6\\
 &+ \Sq^{13}\Sq^5 g_3^3\\
\delta(g_6^{22}) &= g_4^{21}\\
 &+ (\Sq^6\Sq^2 + \Sq^8 + \Sq^7\Sq^1) g_4^{13}\\
 &+ \Sq^{10} g_4^{11}\\
 &+ (\Sq^{13}\Sq^4 + \Sq^{15}\Sq^2 + \Sq^{17}) g_4^4\\
\delta(g_7^{22}) &= (\Sq^{13}\Sq^3 + \Sq^{14}\Sq^2 + \Sq^{16}) g_5^5\\
\ \\
\delta(g_5^{23}) &=\Sq^4 g_3^{18}\\
 &+ \Sq^6\Sq^3\Sq^1 g_3^{12}\\
 &+ (\Sq^{10}\Sq^1 +\Sq^{11}) g_3^{11}\\
 &+ (\Sq^8\Sq^4 + \Sq^9\Sq^3) g_3^{10}\\
 &+ (\Sq^{13}\Sq^3 + \Sq^{15}\Sq^1 + \Sq^9\Sq^4\Sq^2\Sq^1 + \Sq^{11}\Sq^5+ \Sq^{14}\Sq^2 + \Sq^{12}\Sq^4) g_3^6\\
 &+ (\Sq^{16}\Sq^3 + \Sq^{13}\Sq^6 + \Sq^{15}\Sq^4) g_3^3\\
\delta(g_6^{23}) &= g_4^{22}\\
 &+ \Sq^9 g_4^{13}\\
 &+ (\Sq^{10}\Sq^1 + \Sq^{11} + \Sq^8\Sq^3) g_4^{11}\\
 &+ (\Sq^{16}\Sq^2 + \Sq^{14}\Sq^4) g_4^4\\
\delta(g_7^{23}) &= (\Sq^6 + \Sq^4\Sq^2) g_5^{16}\\
 &+ \Sq^7\Sq^1 g_5^{14}\\
 &+ \Sq^{15}\Sq^2 g_5^5\\
\delta(g_8^{23}) &= \Sq^5 g_6^{17}\\
 &+ (\Sq^{14}\Sq^2 +\Sq^{13}\Sq^3) g_6^6\\
\ \\
\delta(g_3^{24}) &=\Sq^{11}\Sq^4 g_1^8\\
 &+ (\Sq^{19} + \Sq^{17}\Sq^2) g_1^4\\
 &+ (\Sq^{16}\Sq^4\Sq^1 + \Sq^{14}\Sq^4\Sq^2\Sq^1 + \Sq^{17}\Sq^4 + \Sq^{21}) g_1^2\\
 &+ (\Sq^{15}\Sq^7 + \Sq^{14}\Sq^6\Sq^2 + \Sq^{18}\Sq^4 + \Sq^{22} + \Sq^{20}\Sq^2 + \Sq^{15}\Sq^5\Sq^2)g_1^1\\
\delta(g_4^{24}) &= \Sq^5 g_2^{18}\\
 &+ (\Sq^{12}\Sq^1 + \Sq^{10}\Sq^2\Sq^1 + \Sq^9\Sq^4) g_2^{10}\\
 &+ (\Sq^{12}\Sq^2 + \Sq^8\Sq^4\Sq^2 + \Sq^{11}\Sq^3) g_2^9\\
 &+ (\Sq^{13}\Sq^2 + \Sq^{14}\Sq^1)g_2^8\\
 &+ (\Sq^{12}\Sq^4\Sq^2 + \Sq^{16}\Sq^2 + \Sq^{12}\Sq^6 + \Sq^{14}\Sq^4 + \Sq^{11}\Sq^5\Sq^2 + \Sq^{15}\Sq^3 + \Sq^{13}\Sq^5) g_2^5\\
 &+ (\Sq^{15}\Sq^4 + \Sq^{19} + \Sq^{14}\Sq^4\Sq^1 + \Sq^{13}\Sq^6 +
 \Sq^{18}\Sq^1 + \Sq^{13}\Sq^4\Sq^2 + \Sq^{15}\Sq^3\Sq^1\\
 &\ \ \ + \Sq^{17}\Sq^2 + \Sq^{12}\Sq^6\Sq^1)g_2^4\\
 &+ (\Sq^{16}\Sq^5 + \Sq^{14}\Sq^5\Sq^2 + \Sq^{18}\Sq^3 + \Sq^{15}\Sq^6 +\Sq^{13}\Sq^6\Sq^2 + \Sq^{15}\Sq^4\Sq^2 + \Sq^{21}) g_2^2
\end{array}
$$

$$
\begin{array}{rl}
\delta(g_7^{24}) &= g_5^{23}\\
 &+ \Sq^4 g_5^{19}\\
 &+ (\Sq^5\Sq^2 + \Sq^7) g_5^{16}\\
 &+ (\Sq^9 + \Sq^8\Sq^1) g_5^{14}\\
 &+ (\Sq^{12}\Sq^6 + \Sq^{16}\Sq^2 + \Sq^{13}\Sq^5 +\Sq^{14}\Sq^4)g_5^5\\
\ \\
\delta(g_5^{25}) &=\Sq^4 g_3^{20}\\
 &+ \Sq^6 g_3^{18}\\
 &+ \Sq^7 g_3^{17}\\
 &+ (\Sq^8\Sq^4 + \Sq^9\Sq^2\Sq^1 + \Sq^9\Sq^3 +\Sq^{10}\Sq^2)g_3^{12}\\
 &+ (\Sq^9\Sq^4 + \Sq^{10}\Sq^2\Sq^1 + \Sq^{12}\Sq^1) g_3^{11}\\ 
 &+ (\Sq^{10}\Sq^5\Sq^2\Sq^1 + \Sq^{12}\Sq^6 + \Sq^{13}\Sq^5 + \Sq^{18} +\Sq^{13}\Sq^4\Sq^1 + \Sq^{16}\Sq^2 + \Sq^{17}\Sq^1\\
 &\ \ \ + \Sq^{14}\Sq^4 + \Sq^{15}\Sq^3) g_3^6\\
 &+ (\Sq^{18}\Sq^3 + \Sq^{19}\Sq^2 + \Sq^{16}\Sq^5) g_3^3\\
\delta(g_8^{25}) &= \Sq^7 g_6^{17}\\
\ \\
\delta(g_4^{26}) &= (\Sq^4\Sq^2\Sq^1 + \Sq^6\Sq^1) g_2^{18}\\
 &+ (\Sq^8 +\Sq^6\Sq^2) g_2^{17}\\
 &+ (\Sq^{15} + \Sq^{14}\Sq^1 + \Sq^8\Sq^4\Sq^2\Sq^1) g_2^{10}\\
 &+ (\Sq^{12}\Sq^4 + \Sq^{14}\Sq^2) g_2^9\\
 &+ (\Sq^{13}\Sq^4 + \Sq^{14}\Sq^3 + \Sq^{10}\Sq^5\Sq^2 + \Sq^{15}\Sq^2)g_2^8\\
 &+ (\Sq^{14}\Sq^6 +\Sq^{18}\Sq^2 + \Sq^{15}\Sq^5 + \Sq^{12}\Sq^6\Sq^2)g_2^5\\
 &+ (\Sq^{18}\Sq^2\Sq^1 + \Sq^{20}\Sq^1 + \Sq^{17}\Sq^4 + \Sq^{19}\Sq^2 + \Sq^{13}\Sq^6\Sq^2 + \Sq^{14}\Sq^7 + \Sq^{15}\Sq^6) g_2^4\\
 &+ (\Sq^{17}\Sq^6 + \Sq^{14}\Sq^7\Sq^2 + \Sq^{19}\Sq^4 + \Sq^{16}\Sq^5\Sq^2 + \Sq^{14}\Sq^6\Sq^3 + \Sq^{15}\Sq^6\Sq^2 + \Sq^{20}\Sq^3\\
 &\ \ \  + \Sq^{17}\Sq^4\Sq^2) g_2^2\\
\delta(g_5^{26}) &=\Sq^5 g_3^{20}\\
 &+ \Sq^5\Sq^2 g_3^{18}\\
 &+ \Sq^6\Sq^2 g_3^{17}\\
 &+ (\Sq^{10}\Sq^2\Sq^1 + \Sq^{10}\Sq^3 +\Sq^8\Sq^4\Sq^1)g_3^{12}\\
 &+ ( \Sq^{13}\Sq^2 + \Sq^{11}\Sq^4 + \Sq^9\Sq^4\Sq^2 + \Sq^{12}\Sq^3 + \Sq^{10}\Sq^5) g_3^{10}\\
 &+ (\Sq^{17}\Sq^2 + \Sq^{11}\Sq^5\Sq^2\Sq^1 + \Sq^{16}\Sq^2\Sq^1 +
 \Sq^{19} + \Sq^{15}\Sq^4 + \Sq^{14}\Sq^4\Sq^1 + \Sq^{12}\Sq^6\Sq^1\\
 &\ \ \ + \Sq^{12}\Sq^4\Sq^2\Sq^1 + \Sq^{12}\Sq^5\Sq^2) g_3^6\\
 &+ (\Sq^{18}\Sq^4 + \Sq^{14}\Sq^6\Sq^2 + \Sq^{16}\Sq^6 + \Sq^{19}\Sq^3 + \Sq^{17}\Sq^5 + \Sq^{13}\Sq^6\Sq^3) g_3^3\\
\delta(g_6^{26}) &= \Sq^3 {}'g_4^{22}\\
 &+ \Sq^3 g_4^{22}\\
 &+ \Sq^4 g_4^{21}\\
 &+ \Sq^6 g_4^{19}\\
 &+ (\Sq^{10}\Sq^2 + \Sq^9\Sq^2\Sq^1 + \Sq^{12}) g_4^{13}\\
 &+(\Sq^{13}\Sq^1 + \Sq^{14} + \Sq^{11}\Sq^3 + \Sq^{12}\Sq^2)g_4^{11}\\
 &+ (\Sq^{17}\Sq^4 + \Sq^{15}\Sq^6 + \Sq^{18}\Sq^3 + \Sq^{21}) g_4^4\\
\delta(g_9^{26}) &=(\Sq^6\Sq^1 + \Sq^4\Sq^2\Sq^1) g_7^{18}\\
 &+(\Sq^{15}\Sq^3 + \Sq^{16}\Sq^2) g_7^7\\
\ \\
\delta(g_4^{27}) &=\Sq^4\Sq^2 g_2^{20}\\
 &+ (\Sq^7\Sq^2+ \Sq^9) g_2^{17}\\
 &+ \Sq^{10}g_2^{16}\\
 &+ (\Sq^{12}\Sq^4 +\Sq^{11}\Sq^4\Sq^1 + \Sq^{13}\Sq^2\Sq^1 + \Sq^{16})g_2^{10}\\
 &+ (\Sq^{17} + \Sq^{10}\Sq^5\Sq^2 + \Sq^{13}\Sq^4 + \Sq^{12}\Sq^5 +\Sq^{15}\Sq^2) g_2^9\\
 &+ (\Sq^{12}\Sq^4\Sq^2 + \Sq^{14}\Sq^4) g_2^8\\
 &+ (\Sq^{15}\Sq^6 + \Sq^{19}\Sq^2 + \Sq^{12}\Sq^6\Sq^3 + \Sq^{13}\Sq^6\Sq^2) g_2^5\\
 &+ (\Sq^{17}\Sq^4\Sq^1 + \Sq^{12}\Sq^6\Sq^3\Sq^1 + \Sq^{20}\Sq^2 +\Sq^{15}\Sq^7 + \Sq^{14}\Sq^7\Sq^1) g_2^4\\
 &+ (\Sq^{15}\Sq^6\Sq^3 +\Sq^{18}\Sq^6 + \Sq^{16}\Sq^8 + \Sq^{20}\Sq^4 + \Sq^{18}\Sq^4\Sq^2 +\Sq^{15}\Sq^7\Sq^2) g_2^2
\end{array}
$$

\ 

$$
\begin{array}{rl}
\delta(g_5^{28}) &= (\Sq^6\Sq^1 + \Sq^7) g_3^{20}\\
 &+ \Sq^9 g_3^{18}\\
 &+ \Sq^7\Sq^3g_3^{17}\\
 &+ (\Sq^{12}\Sq^2\Sq^1 + \Sq^9\Sq^4\Sq^2 +\Sq^8\Sq^4\Sq^2\Sq^1 + \Sq^{11}\Sq^3\Sq^1 + \Sq^{15} + \Sq^{14}\Sq^1 +\Sq^{12}\Sq^3) g_3^{12}\\
 &+ \Sq^{12}\Sq^4 g_3^{11}\\
 &+ (\Sq^{14}\Sq^3 + \Sq^{11}\Sq^4\Sq^2) g_3^{10}\\
 &+ ( \Sq^{14}\Sq^4\Sq^2\Sq^1 + \Sq^{15}\Sq^6 + \Sq^{21} +\Sq^{13}\Sq^6\Sq^2 + \Sq^{19}\Sq^2 + \Sq^{18}\Sq^3 + \Sq^{15}\Sq^4\Sq^2\\
 &\ \ \ + \Sq^{17}\Sq^4 + \Sq^{14}\Sq^5\Sq^2 + \Sq^{18}\Sq^2\Sq^1 + \Sq^{12}\Sq^6\Sq^3 + \Sq^{14}\Sq^7 + \Sq^{12}\Sq^6\Sq^2\Sq^1) g_3^6\\
 &+ (\Sq^{20}\Sq^4 + \Sq^{24} + \Sq^{18}\Sq^6 + \Sq^{18}\Sq^4\Sq^2 + \Sq^{16}\Sq^8)g_3^3\\
\delta(g_9^{28}) &=\Sq^4 g_7^{23}\\
 &+ (\Sq^{20} + \Sq^{18}\Sq^2) g_7^7\\
\delta(g_{10}^{28}) &= 0\\
\ \\
\delta(g_5^{29}) &=\Sq^4\Sq^2 g_3^{22}\\
 &+ (\Sq^7\Sq^1 + \Sq^6\Sq^2)g_3^{20}\\
 &+ \Sq^{10} g_3^{18}\\
 &+ (\Sq^{11} + \Sq^9\Sq^2) g_3^{17}\\
 &+ (\Sq^{12}\Sq^4 + \Sq^9\Sq^4\Sq^2\Sq^1 + \Sq^{16} + \Sq^{10}\Sq^4\Sq^2 + \Sq^{12}\Sq^3\Sq^1 + \Sq^{15}\Sq^1 + \Sq^{14}\Sq^2)g_3^{12}\\
 &+ (\Sq^{17} + \Sq^{16}\Sq^1) g_3^{11}\\
 &+ (\Sq^{11}\Sq^5\Sq^2 + \Sq^{12}\Sq^4\Sq^2 + \Sq^{13}\Sq^5 +\Sq^{18})g_3^{10}\\
 &+ (\Sq^{19}\Sq^2\Sq^1 + \Sq^{15}\Sq^6\Sq^1 + \Sq^{15}\Sq^4\Sq^2\Sq^1 + \Sq^{14}\Sq^5\Sq^2\Sq^1 + \Sq^{12}\Sq^6\Sq^3\Sq^1\\
 &\ \ \ + \Sq^{15}\Sq^5\Sq^2 + \Sq^{17}\Sq^5 + \Sq^{19}\Sq^3 + \Sq^{22}) g_3^6\\
 &+ ( \Sq^{16}\Sq^6\Sq^3 + \Sq^{20}\Sq^5 + \Sq^{17}\Sq^6\Sq^2 + \Sq^{15}\Sq^7\Sq^3 + \Sq^{16}\Sq^7\Sq^2 + \Sq^{18}\Sq^5\Sq^2\\
 &\ \ \ + \Sq^{19}\Sq^6) g_3^3\\
\delta(g_6^{29}) &= (\Sq^{12}\Sq^2\Sq^1 + \Sq^{11}\Sq^4 + \Sq^{13}\Sq^2 + \Sq^{14}\Sq^1 + \Sq^{15} + \Sq^{11}\Sq^3\Sq^1) g_4^{13}\\
 &+ (\Sq^{12}\Sq^5 + \Sq^{15}\Sq^2 + \Sq^{13}\Sq^4 + \Sq^{16}\Sq^1 + \Sq^{17} + \Sq^{11}\Sq^4\Sq^2 + \Sq^{14}\Sq^3)g_4^{11}\\ 
 &+ ( \Sq^{17}\Sq^5\Sq^2 + \Sq^{22}\Sq^2 + \Sq^{15}\Sq^7\Sq^2 + \Sq^{18}\Sq^4\Sq^2 + \Sq^{20}\Sq^4) g_4^4\\
\delta(g_{10}^{29}) &= (\Sq^5 + \Sq^4\Sq^1) g_8^{23}\\
 &+ (\Sq^{18}\Sq^2 + \Sq^{20}) g_8^8\\
\ \\
\delta(g_7^{30}) &= \Sq^3 g_5^{26}\\
 &+ \Sq^4\Sq^2 g_5^{23}\\
 &+ \Sq^7 g_5^{22}\\
 &+ \Sq^9 g_5^{20}\\
 &+ \Sq^9{}'g_5^{20}\\
 &+ \Sq^8\Sq^2 g_5^{19}\\
 &+ (\Sq^{10}\Sq^3 + \Sq^{11}\Sq^2) g_5^{16}\\
 &+ (\Sq^{15} + \Sq^{12}\Sq^2\Sq^1 + \Sq^{10}\Sq^5 + \Sq^8\Sq^4\Sq^2\Sq^1) g_5^{14}\\
 &+ (\Sq^{21}\Sq^3 + \Sq^{19}\Sq^5 + \Sq^{18}\Sq^6 + \Sq^{17}\Sq^7 +\Sq^{22}\Sq^2 + \Sq^{24})g_5^5\\
\delta(g_8^{30}) &= \Sq^2\Sq^1 g_6^{26}\\
 &+ \Sq^6 g_6^{23}\\
 &+ (\Sq^6\Sq^1 + \Sq^4\Sq^2\Sq^1) g_6^{22}\\
 &+ (\Sq^7\Sq^2 + \Sq^8\Sq^1 + \Sq^6\Sq^3)g_6^{20}\\
 &+ (\Sq^{10}\Sq^2 + \Sq^9\Sq^3 + \Sq^8\Sq^4) g_6^{17}\\
 &+ (\Sq^{13} + \Sq^{12}\Sq^1 + \Sq^{10}\Sq^2\Sq^1) g_6^{16}\\
 &+ ( \Sq^{21}\Sq^2 + \Sq^{17}\Sq^6 + \Sq^{18}\Sq^5 + \Sq^{16}\Sq^7 +\Sq^{17}\Sq^4\Sq^2) g_6^6\\
\delta(g_{11}^{30}) &= \Sq^2\Sq^1 g_9^{26}
\end{array}
$$

\ 

$$
\begin{array}{rl}
\delta(g_8^{31}) &= \Sq^7 g_6^{23}\\
 &+ (\Sq^5\Sq^2\Sq^1 + \Sq^7\Sq^1) g_6^{22}\\
 &+ (\Sq^7\Sq^3 + \Sq^{10} + \Sq^8\Sq^2) g_6^{20}\\
 &+ (\Sq^{10}\Sq^3 + \Sq^{13}) g_6^{17}\\
 &+ (\Sq^8\Sq^4\Sq^2 + \Sq^{13}\Sq^1 + \Sq^{14} +\Sq^{11}\Sq^2\Sq^1)g_6^{16}\\
 &+ ( \Sq^{17}\Sq^5\Sq^2 + \Sq^{22}\Sq^2 + \Sq^{21}\Sq^3 + \Sq^{24} +\Sq^{18}\Sq^6 + \Sq^{17}\Sq^7) g_6^6\\
\delta(g_9^{31}) &= g_7^{30}\\
 &+ \Sq^4\Sq^2 g_7^{24}\\
 &+ (\Sq^6\Sq^1 + \Sq^4\Sq^2\Sq^1 + \Sq^7) g_7^{23}\\
 &+ (\Sq^{12} + \Sq^8\Sq^4 + \Sq^9\Sq^2\Sq^1) g_7^{18}\\
 &+ (\Sq^{20}\Sq^3 + \Sq^{18}\Sq^5) g_7^7\\
\ \\
\delta(g_2^{32}) &= 0\\
\delta(g_6^{32}) &= \Sq^4 g_4^{27}\\
 &+ \Sq^5 g_4^{26}\\
 &+ (\Sq^8\Sq^1 + \Sq^9 + \Sq^6\Sq^3) {}'g_4^{22}\\
 &+ \Sq^7\Sq^2 g_4^{22}\\
 &+ (\Sq^{10} + \Sq^8\Sq^2) g_4^{21}\\
 &+ (\Sq^{12} + \Sq^8\Sq^4) g_4^{19}\\
 &+ \Sq^{13} g_4^{18}\\
 &+ (\Sq^{16}\Sq^2 + \Sq^{17}\Sq^1 + \Sq^{13}\Sq^5 + \Sq^{12}\Sq^4\Sq^2 +\Sq^{15}\Sq^2\Sq^1 + \Sq^{14}\Sq^3\Sq^1) g_4^{13}\\
 &+ (\Sq^{15}\Sq^4\Sq^1 + \Sq^{14}\Sq^6 + \Sq^{16}\Sq^4 + \Sq^{14}\Sq^4\Sq^2 + \Sq^{12}\Sq^6\Sq^2 + \Sq^{13}\Sq^6\Sq^1)g_4^{11}\\
 &+ (\Sq^{24}\Sq^3 + \Sq^{17}\Sq^8\Sq^2 + \Sq^{20}\Sq^5\Sq^2 +\Sq^{19}\Sq^6\Sq^2 + \Sq^{18}\Sq^9 + \Sq^{27}) g_4^4\\
\delta(g_9^{32}) &= \Sq^7 g_7^{24}\\
 &+ \Sq^8 g_7^{23}\\
 &+ (\Sq^{11}\Sq^2 + \Sq^{12}\Sq^1 + \Sq^8\Sq^4\Sq^1 +\Sq^9\Sq^4)g_7^{18}\\
 &+ (\Sq^{18}\Sq^4\Sq^2 + \Sq^{21}\Sq^3 + \Sq^{22}\Sq^2 +\Sq^{24} + \Sq^{20}\Sq^4) g_7^7\\
\delta({}'g_9^{32}) &= (\Sq^7 + \Sq^5\Sq^2) g_7^{24}\\
 &+ \Sq^8 g_7^{23}\\
 &+ (\Sq^8\Sq^4\Sq^1 + \Sq^{10}\Sq^2\Sq^1) g_7^{18}\\
 &+ (\Sq^{18}\Sq^4\Sq^2 + \Sq^{19}\Sq^5 + \Sq^{22}\Sq^2+ \Sq^{20}\Sq^4)g_7^7\\
\delta(g_{10}^{32}) &= g_8^{31}\\
 &+ \Sq^6 g_8^{25}\\
 &+ (\Sq^7\Sq^1 + \Sq^6\Sq^2) g_8^{23}\\
 &+ (\Sq^{21}\Sq^2 + \Sq^{19}\Sq^4 + \Sq^{23}) g_8^8\\
\ \\
\delta(g_2^{33}) &= 0\\
\delta(g_3^{33}) &= g_1^{32}\\
 &+ \Sq^{24} g_1^8\\
 &+ (\Sq^{28} + \Sq^{25}\Sq^3) g_1^4\\
 &+ (\Sq^{29}\Sq^1 + \Sq^{30} + \Sq^{23}\Sq^7 +\Sq^{22}\Sq^7\Sq^1 + \Sq^{25}\Sq^5 + \Sq^{23}\Sq^6\Sq^1 +\Sq^{23}\Sq^4\Sq^2\Sq^1) g_1^2\\
 &+ (\Sq^{29}\Sq^2 + \Sq^{20}\Sq^8\Sq^3 + \Sq^{22}\Sq^7\Sq^2 + \Sq^{27}\Sq^4 + \Sq^{21}\Sq^7\Sq^3 + \Sq^{21}\Sq^8\Sq^2 + \Sq^{26}\Sq^5\\
 &\ \ \ + \Sq^{28}\Sq^3 + \Sq^{19}\Sq^9\Sq^3 + \Sq^{25}\Sq^6 + \Sq^{19}\Sq^8\Sq^4) g_1^1\\
\delta(g_7^{33}) &= \Sq^4\Sq^2 g_5^{26}\\
 &+ \Sq^7 g_5^{25}\\
 &+ (\Sq^7\Sq^2 + \Sq^8\Sq^1 + \Sq^6\Sq^3 + \Sq^6\Sq^2\Sq^1)g_5^{23}\\
 &+ \Sq^8\Sq^2 g_5^{22}\\
 &+ (\Sq^8\Sq^4 + \Sq^{11}\Sq^1) g_5^{20}\\
 &+ \Sq^{10}\Sq^2 {}'g_5^{20}\\
 &+ (\Sq^{13}\Sq^5 + \Sq^{15}\Sq^2\Sq^1) g_5^{14}\\
 &+ (\Sq^{18}\Sq^6\Sq^3 + \Sq^{18}\Sq^7\Sq^2 + \Sq^{17}\Sq^8\Sq^2 + \Sq^{22}\Sq^5 + \Sq^{23}\Sq^4 + \Sq^{19}\Sq^6\Sq^2\\
 &\ \ \ + \Sq^{18}\Sq^9)g_5^5
\end{array}
$$

$$
\begin{array}{rl}
\delta(g_8^{33}) &= \Sq^2\Sq^1 g_6^{29}\\
 &+ \Sq^6 g_6^{26}\\
 &+ (\Sq^7\Sq^2 + \Sq^9 + \Sq^6\Sq^3) g_6^{23}\\
 &+ (\Sq^{10} + \Sq^8\Sq^2) g_6^{22}\\
 &+ \Sq^{11} g_6^{21}\\
 &+ (\Sq^8\Sq^4 + \Sq^{10}\Sq^2) g_6^{20}\\
 &+ (\Sq^{13}\Sq^2 + \Sq^{11}\Sq^4 + \Sq^{10}\Sq^5 + \Sq^9\Sq^4\Sq^2 + \Sq^{15}) g_6^{17}\\
 &+ (\Sq^{15}\Sq^1 + \Sq^{13}\Sq^2\Sq^1 + \Sq^9\Sq^4\Sq^2\Sq^1)g_6^{16}\\
 &+ (\Sq^{20}\Sq^6 + \Sq^{18}\Sq^8 + \Sq^{19}\Sq^5\Sq^2 + \Sq^{21}\Sq^5 + \Sq^{23}\Sq^3) g_6^6\\
\delta(g_{10}^{33}) &= (\Sq^7 + \Sq^6\Sq^1) g_8^{25}\\
 &+ (\Sq^9 + \Sq^8\Sq^1) g_8^{23}\\
 &+ \Sq^{19}\Sq^5 g_8^8\\
\ \\
\delta(g_2^{34}) &= 0\\
\delta(g_3^{34}) &= (\Sq^{21}\Sq^8 + \Sq^{22}\Sq^5\Sq^2) g_1^4\\
 &+ (\Sq^{19}\Sq^8\Sq^4 + \Sq^{23}\Sq^8 + \Sq^{16}\Sq^8\Sq^4\Sq^2\Sq^1 + \Sq^{21}\Sq^{10} + \Sq^{24}\Sq^4\Sq^2\Sq^1 + \Sq^{30}\Sq^1\\
 &\ \ \ + \Sq^{21}\Sq^7\Sq^2\Sq^1 + \Sq^{20}\Sq^{10}\Sq^1 + \Sq^{23}\Sq^7\Sq^1) g_1^2\\
 &+ (\Sq^{25}\Sq^7 + \Sq^{22}\Sq^8\Sq^2 + \Sq^{25}\Sq^5\Sq^2 + \Sq^{30}\Sq^2 + \Sq^{20}\Sq^9\Sq^3 + \Sq^{29}\Sq^3 + \Sq^{23}\Sq^6\Sq^3\\
 &\ \ \ + \Sq^{23}\Sq^9 + \Sq^{20}\Sq^8\Sq^4 + \Sq^{21}\Sq^8\Sq^3 + \Sq^{26}\Sq^4\Sq^2 + \Sq^{18}\Sq^8\Sq^4\Sq^2 + \Sq^{20}\Sq^{10}\Sq^2\\
 &\ \ \ + \Sq^{24}\Sq^8 + \Sq^{32}) g_1^1\\
\delta({}'g_3^{34}) &= (\Sq^{26}\Sq^3 + \Sq^{21}\Sq^8 + \Sq^{27}\Sq^2 + \Sq^{29}) g_1^4\\
 &+ (\Sq^{23}\Sq^8 + \Sq^{28}\Sq^2\Sq^1 + \Sq^{21}\Sq^7\Sq^2\Sq^1 + \Sq^{22}\Sq^6\Sq^2\Sq^1 + \Sq^{27}\Sq^4 + \Sq^{24}\Sq^7\\
 &\ \ \ + \Sq^{19}\Sq^8\Sq^4 + \Sq^{16}\Sq^8\Sq^4\Sq^2\Sq^1) g_1^2\\
 &+ (\Sq^{24}\Sq^6\Sq^2 + \Sq^{30}\Sq^2 + \Sq^{27}\Sq^5 + \Sq^{24}\Sq^8 + \Sq^{21}\Sq^8\Sq^3 + \Sq^{25}\Sq^5\Sq^2 + \Sq^{23}\Sq^6\Sq^3\\
 &\ \ \ + \Sq^{18}\Sq^8\Sq^4\Sq^2 + \Sq^{26}\Sq^4\Sq^2 + \Sq^{22}\Sq^8\Sq^2 + \Sq^{20}\Sq^8\Sq^4 + \Sq^{32} + \Sq^{29}\Sq^3\\
 & + \Sq^{23}\Sq^9) g_1^1\\
\delta(g_4^{34}) &= g_2^{33}\\
 &+ \Sq^{13} g_2^{20}\\
 &+ \Sq^{15} g_2^{18}\\
 &+ (\Sq^{16} + \Sq^{11}\Sq^5 + \Sq^{10}\Sq^4\Sq^2) g_2^{17}\\
 &+ \Sq^{12}\Sq^5 g_2^{16}\\
 &+ (\Sq^{23} + \Sq^{17}\Sq^6 + \Sq^{16}\Sq^4\Sq^2\Sq^1 + \Sq^{15}\Sq^7\Sq^1 + \Sq^{16}\Sq^6\Sq^1 + \Sq^{18}\Sq^4\Sq^1) g_2^{10}\\
 &+ (\Sq^{21}\Sq^3 + \Sq^{18}\Sq^6 + \Sq^{16}\Sq^6\Sq^2 + \Sq^{15}\Sq^7\Sq^2 + \Sq^{17}\Sq^5\Sq^2 + \Sq^{15}\Sq^6\Sq^3\\
 & + \Sq^{19}\Sq^5) g_2^9\\
 &+ (\Sq^{25} + \Sq^{23}\Sq^2 + \Sq^{21}\Sq^4 + \Sq^{19}\Sq^4\Sq^2 + \Sq^{18}\Sq^5\Sq^2 + \Sq^{22}\Sq^3) g_2^8\\
 &+ (\Sq^{26}\Sq^2 + \Sq^{20}\Sq^6\Sq^2 + \Sq^{22}\Sq^4\Sq^2 + \Sq^{21}\Sq^5\Sq^2 + \Sq^{20}\Sq^8 + \Sq^{19}\Sq^9 + \Sq^{25}\Sq^3\\
 &\ \ \ + \Sq^{24}\Sq^4 + \Sq^{17}\Sq^8\Sq^3 + \Sq^{16}\Sq^8\Sq^4 + \Sq^{18}\Sq^8\Sq^2)g_2^5\\
 &+ (\Sq^{19}\Sq^6\Sq^3\Sq^1 + \Sq^{20}\Sq^9 + \Sq^{27}\Sq^2 + \Sq^{20}\Sq^7\Sq^2 + \Sq^{23}\Sq^4\Sq^2 + \Sq^{28}\Sq^1 + \Sq^{29}\\
 &\ \ \ + \Sq^{20}\Sq^8\Sq^1 + \Sq^{22}\Sq^7 + \Sq^{22}\Sq^6\Sq^1 + \Sq^{24}\Sq^4\Sq^1 + \Sq^{22}\Sq^4\Sq^2\Sq^1\\
 &\ \ \ + \Sq^{18}\Sq^8\Sq^2\Sq^1 + \Sq^{25}\Sq^4 + \Sq^{19}\Sq^7\Sq^2\Sq^1 + \Sq^{21}\Sq^5\Sq^2\Sq^1 + \Sq^{26}\Sq^2\Sq^1) g_2^4\\
 &+ ( \Sq^{20}\Sq^9\Sq^2 + \Sq^{19}\Sq^8\Sq^4 + \Sq^{21}\Sq^7\Sq^3 + \Sq^{26}\Sq^5 + \Sq^{21}\Sq^{10} + \Sq^{31} + \Sq^{20}\Sq^8\Sq^3\\
 &\ \ \ + \Sq^{18}\Sq^9\Sq^4 + \Sq^{19}\Sq^9\Sq^3 + \Sq^{21}\Sq^8\Sq^2 + \Sq^{22}\Sq^7\Sq^2) g_2^2\\
\delta(g_8^{34}) &=(\Sq^6\Sq^1 + \Sq^4\Sq^2\Sq^1 + \Sq^7) g_6^{26}\\
 &+ (\Sq^9\Sq^2 + \Sq^{10}\Sq^1 + \Sq^8\Sq^2\Sq^1) g_6^{22}\\
 &+ \Sq^{12} g_6^{21}\\
 &+ (\Sq^{11}\Sq^5 + \Sq^{13}\Sq^3) g_6^{17}\\
 &+ (\Sq^{13}\Sq^4 + \Sq^{12}\Sq^4\Sq^1 + \Sq^{15}\Sq^2) g_6^{16}\\
 &+ (\Sq^{18}\Sq^9 + \Sq^{21}\Sq^6 + \Sq^{22}\Sq^5 + \Sq^{27} + \Sq^{19}\Sq^6\Sq^2 + \Sq^{24}\Sq^3)g_6^6\\
\delta(g_9^{34}) &=g_7^{33}\\
 &+ (\Sq^6\Sq^3 + \Sq^9 + \Sq^6\Sq^2\Sq^1) g_7^{24}\\
 &+ (\Sq^{14}\Sq^1 + \Sq^{13}\Sq^2 + \Sq^{15} + \Sq^{11}\Sq^4 + \Sq^8\Sq^4\Sq^2\Sq^1 + \Sq^{11}\Sq^3\Sq^1) g_7^{18}\\
 &+ (\Sq^{19}\Sq^5\Sq^2 + \Sq^{20}\Sq^4\Sq^2 + \Sq^{24}\Sq^2 + \Sq^{18}\Sq^8 + \Sq^{21}\Sq^5 + \Sq^{23}\Sq^3 + \Sq^{26}) g_7^7
\end{array}
$$

$$
\begin{array}{rl}
\delta(g_{10}^{34}) &= g_8^{33}\\
 &+ (\Sq^6\Sq^2 + \Sq^8 + \Sq^5\Sq^2\Sq^1 + \Sq^7\Sq^1)g_8^{25}\\
 &+ \Sq^8\Sq^2 g_8^{23}\\
 &+ (\Sq^{19}\Sq^6 + \Sq^{21}\Sq^4 + \Sq^{25} + \Sq^{20}\Sq^5) g_8^8\\
\delta(g_{11}^{34}) &= \Sq^1 {}'g_9^{32}\\
 &+ (\Sq^{21}\Sq^3 + \Sq^{22}\Sq^2 + \Sq^{24}) g_9^9\\
\ \\
\delta(g_4^{35}) &= \Sq^{17} g_2^{17}\\
 &+ (\Sq^{18} + \Sq^{16}\Sq^2) g_2^{16}\\
 &+ (\Sq^{18}\Sq^6 + \Sq^{16}\Sq^7\Sq^1) g_2^{10}\\
 &+ \Sq^{19}\Sq^6 g_2^9\\
 &+ (\Sq^{18}\Sq^8 + \Sq^{19}\Sq^5\Sq^2 + \Sq^{24}\Sq^2 + \Sq^{18}\Sq^7\Sq^1 + \Sq^{16}\Sq^8\Sq^2 + \Sq^{26}) g_2^8\\
 &+ (\Sq^{25}\Sq^4 + \Sq^{21}\Sq^6\Sq^2 + \Sq^{24}\Sq^5 + \Sq^{26}\Sq^3 + \Sq^{21}\Sq^8 + \Sq^{23}\Sq^6 + \Sq^{18}\Sq^9\Sq^2\\
 &\ \ \ + \Sq^{19}\Sq^7\Sq^3 + \Sq^{20}\Sq^7\Sq^2) g_2^5\\
 &+ (\Sq^{23}\Sq^7 + \Sq^{20}\Sq^8\Sq^2 + \Sq^{20}\Sq^7\Sq^2\Sq^1 + \Sq^{23}\Sq^4\Sq^2\Sq^1 + \Sq^{27}\Sq^2\Sq^1\\
 &\ \ \ + \Sq^{16}\Sq^8\Sq^4\Sq^2 + \Sq^{22}\Sq^6\Sq^2 + \Sq^{23}\Sq^6\Sq^1 + \Sq^{19}\Sq^8\Sq^2\Sq^1 + \Sq^{28}\Sq^2 + \Sq^{25}\Sq^5\\
 &\ \ \ + \Sq^{21}\Sq^7\Sq^2) g_2^4\\
 &+ (\Sq^{30}\Sq^2 + \Sq^{20}\Sq^9\Sq^3 + \Sq^{22}\Sq^8\Sq^2 + \Sq^{23}\Sq^6\Sq^3 + \Sq^{25}\Sq^5\Sq^2 + \Sq^{24}\Sq^8\\
 &\ \ \ + \Sq^{26}\Sq^4\Sq^2 + \Sq^{24}\Sq^6\Sq^2 + \Sq^{21}\Sq^9\Sq^2 + \Sq^{23}\Sq^9 + \Sq^{18}\Sq^8\Sq^4\Sq^2)g_2^2\\
\delta(g_5^{35}) &={}'g_3^{34}\\
 &+ \Sq^9\Sq^4 g_3^{21}\\
 &+ (\Sq^8\Sq^4\Sq^2 + \Sq^{10}\Sq^3\Sq^1 + \Sq^{11}\Sq^2\Sq^1)g_3^{20}\\
 &+ (\Sq^{14}\Sq^2 + \Sq^{12}\Sq^4 + \Sq^{13}\Sq^3) g_3^{18}\\
 &+ (\Sq^{12}\Sq^5 + \Sq^{17} + \Sq^{14}\Sq^2\Sq^1) g_3^{17}\\
 &+ (\Sq^{16}\Sq^4\Sq^2 + \Sq^{14}\Sq^6\Sq^2 + \Sq^{12}\Sq^6\Sq^3\Sq^1 + \Sq^{17}\Sq^5 + \Sq^{14}\Sq^7\Sq^1 + \Sq^{21}\Sq^1\\
 &\ \ \ + \Sq^{19}\Sq^3 + \Sq^{15}\Sq^7) g_3^{12}\\
 &+ ( \Sq^{20}\Sq^2\Sq^1 + \Sq^{17}\Sq^6 + \Sq^{15}\Sq^7\Sq^1 + \Sq^{18}\Sq^5 + \Sq^{22}\Sq^1) g_3^{11}\\
 &+ ( \Sq^{24} + \Sq^{20}\Sq^4 + \Sq^{19}\Sq^5 + \Sq^{17}\Sq^5\Sq^2 + \Sq^{18}\Sq^6) g_3^{10}\\
 &+ (\Sq^{20}\Sq^8 + \Sq^{23}\Sq^4\Sq^1 + \Sq^{24}\Sq^3\Sq^1 + \Sq^{24}\Sq^4 + \Sq^{18}\Sq^6\Sq^3\Sq^1 + \Sq^{21}\Sq^4\Sq^2\Sq^1\\
 &\ \ \ + \Sq^{19}\Sq^8\Sq^1 + \Sq^{21}\Sq^7 + \Sq^{25}\Sq^3 + \Sq^{18}\Sq^7\Sq^3 + \Sq^{19}\Sq^9 + \Sq^{27}\Sq^1 + \Sq^{16}\Sq^8\Sq^4\\
 &\ \ \ + \Sq^{23}\Sq^5 + \Sq^{17}\Sq^8\Sq^3 + \Sq^{21}\Sq^6\Sq^1) g_3^6\\
 &+ ( \Sq^{31} + \Sq^{26}\Sq^5 + \Sq^{25}\Sq^6 + \Sq^{19}\Sq^8\Sq^4 + \Sq^{21}\Sq^7\Sq^3 + \Sq^{19}\Sq^9\Sq^3 + \Sq^{29}\Sq^2\\
 &\ \ \ + \Sq^{24}\Sq^5\Sq^2 + \Sq^{20}\Sq^8\Sq^3 + \Sq^{25}\Sq^4\Sq^2 + \Sq^{22}\Sq^6\Sq^3 + \Sq^{20}\Sq^9\Sq^2 + \Sq^{23}\Sq^8\\
 &\ \ \ + \Sq^{27}\Sq^4 + \Sq^{24}\Sq^7 + \Sq^{28}\Sq^3) g_3^3\\
\delta(g_9^{35}) &=\Sq^4 g_7^{30}\\
 &+ (\Sq^7\Sq^3 + \Sq^6\Sq^3\Sq^1 + \Sq^{10} + \Sq^8\Sq^2)g_7^{24}\\
 &+ (\Sq^{18}\Sq^9 + \Sq^{19}\Sq^6\Sq^2 + \Sq^{23}\Sq^4 + \Sq^{24}\Sq^3 + \Sq^{21}\Sq^6 + \Sq^{19}\Sq^8) g_7^7\\
 &+ (\Sq^{11} + \Sq^8\Sq^2\Sq^1 + \Sq^{10}\Sq^1) g_7^{23}\\
 &+ \Sq^{12} g_7^{22}\\
 &+ (\Sq^{13}\Sq^2\Sq^1 + \Sq^{12}\Sq^4 + \Sq^{12}\Sq^3\Sq^1 + \Sq^{10}\Sq^4\Sq^2 + \Sq^{11}\Sq^4\Sq^1 + \Sq^9\Sq^4\Sq^2\Sq^1) g_7^{18}\\
\delta(g_{10}^{35}) &= g_8^{34}\\
 &+ \Sq^3 g_8^{31}\\
 &+ (\Sq^{10}\Sq^1 + \Sq^{11} + \Sq^9\Sq^2) g_8^{23}\\
 &+ (\Sq^{24}\Sq^2 + \Sq^{19}\Sq^7) g_8^8\\
\delta(g_{11}^{35}) &= \Sq^2 {}'g_9^{32}\\
 &+ \Sq^3 g_9^{31}\\
 &+ \Sq^6 g_9^{28}\\
 &+ (\Sq^8 + \Sq^7\Sq^1) g_9^{26}\\
 &+ \Sq^{25} g_9^9\\
\delta(g_{12}^{35}) &=\Sq^5 g_{10}^{29}\\
 &+ (\Sq^{22}\Sq^2 + \Sq^{21}\Sq^3)g_{10}^{10}\\
\end{array}
$$

\end{document}